\documentclass[12pt]{article}
\usepackage{amsfonts, amssymb, latexsym, graphicx}
\usepackage[all]{xy}
\input xy
\xyoption{all}
\newtheorem{theorem}{Theorem}
\newtheorem{proposition}{Proposition}

\begin{document}
\title{Rational pencils of cubics and configurations of six or seven points in $\mathbb{R}P^2$}
\author{S\'everine Fiedler-Le Touz\'e}
\maketitle
\begin{abstract}
Let six points $1, \dots 6$ lie in general position in the real projective plane and consider the pencil of nodal cubics based at these points, with node at one of them, say $1$. This pencil has five reducible cubics. We call combinatorial cubic a topological type (cubic, points) and combinatorial pencil the cyclic sequence of five combinatorial reducible cubics. Up to the action of the symmetric group $\mathcal{S}_5$ on $\{ 2, \dots 6\}$, there are seven possible combinatorial pencils with node at $1$. Consider now the set of six pencils obtained, making the node to be $1, \dots 6$. Up to the action of $\mathcal{S}_6$ on $\{1, \dots 6\}$, there are four possible lists of six combinatorial pencils. Let seven points $1, \dots 7$ lie in general position in the plane. Up to the action of $\mathcal{S}_7$ on $\{1, \dots 7\}$, there are fourteen possible lists of seven nodal combinatorial cubics passing through the seven points, with respective nodes at $1, \dots 7$.
\end{abstract}

\section{Introduction}

Let $1, \dots, 6$ be six points in the real projective plane and consider the pencil of nodal cubics determined by these points with node at one of them, say $1$.  This pencil has five reducible cubics plus, possibly, some cuspidal cubics. We ignore the latter and consider the sequence of reducible cubics. Let us call combinatorial cubic a topological type (cubic, points), up to the following identification: a loop passing through no other point than the node will be assimilated to an isolated node. Let us
call combinatorial pencil the cyclic sequence of five combinatorial reducible cubics, they are of the form $1m \cup 1ijkl$, where $\{i, j, k, l, m \} = \{ 2, 3, 4, 5, 6 \}$. 
A set of six or seven points is in generic position if no three are aligned and no six are coconic. Consider first six points, five of them determine $10$ lines and one conic, dividing $\mathbb{R}P^2$ in $36$ zones. The {\em (ordered) configuration\/} of $(1, \dots 6)$ is the topological type of the sextuple of points with respect to the lines and conics, described by any one of the six equivalent pieces of information: cyclic ordering of five chosen points on their conic and zone containing the sixth one. Another equivalent information is the list of six combinatorial pencils of cubics based at $1, \dots 6$, with respective nodes at $1, \dots 6$. An {\em unordered configuration\/} of six points is an equivalence class of ordered configurations for the action of the group $\mathcal{S}_6$. 
Let us consider now seven points $1, \dots 7$. For a point $n$ among the seven, we denote by $\hat n$ the ordered configuration of the other six points. An ordered configuration of seven points is the set $\hat 1, \dots \hat 7$ of ordered configurations realized by six of the points. An unordered configuration of seven points is an equivalence class of ordered configurations for the action of $\mathcal{S}_7$. 

\begin{theorem}
Let $1, \dots, 6$ be six generic points in $\mathbb{R}P^2$. Up to the action of $\mathcal{S}_5$ on $\{2, \dots 6\}$, there are seven possible combinatorial pencils of nodal cubics based at these six points with node at $1$.

Up to the action of $\mathcal{S}_6$ on $\{1, \dots 6\}$, there are four possible lists of six combinatorial pencils of nodal cubics based at these six points with respective nodes at $1, \dots 6$. Otherwise stated, six generic points may realize four different unordered configurations. 
\end{theorem}

\begin{theorem}
Seven generic points $1, \dots 7$ in  $\mathbb{R}P^2$ may realize fourteen different unordered configurations. Up to the action of $\mathcal{S}_7$, there are correspondingly fourteen possible lists of seven combinatorial nodal cubics through the seven points, with respective nodes at $1, \dots 7$.
\end{theorem} 

Theorem~1 is proved in section 2.1, Theorem~2 is proved in sections 3.1 and 3.2.
The $14$ configurations have been obtained independently by Arzu Zabun. 
The recent paper \cite{fz} by Finashin and Zabun gives the classification of the Aronhold sets of seven bitangents to real plane quartics with four ovals. There are $14$ of them, corresponding via the del Pezzo surfaces of degree $2$ to the $14$ configurations.

{\sc Acknowledgement:} I am grateful to Arzu Zabun and Sergey Finashin
for pointing out to me that one configuration of Theorem~2 was missing in a previous version (v2) of this paper. 

\section{Configurations of six points}
\subsection{Rational pencils of cubics}

\begin{figure}
\centering
\includegraphics{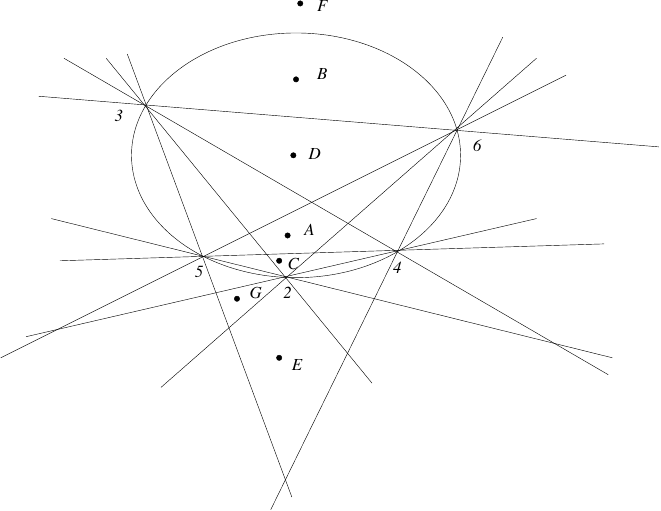}
\caption{\label{conf} The $36$ zones}  
\end{figure}

{\em Proof of Theorem~1:\/}
Up to the action of $\mathcal{S}_5$ on $\{2, \dots 6\}$, we may assume that the six points are disposed as shown in Figure~\ref{conf}, where $2, 5, 3, 6, 4$ lie in this ordering on a conic, and $1$ is in one of the zones $A, \dots G$. 

\begin{figure}
\centering
\includegraphics{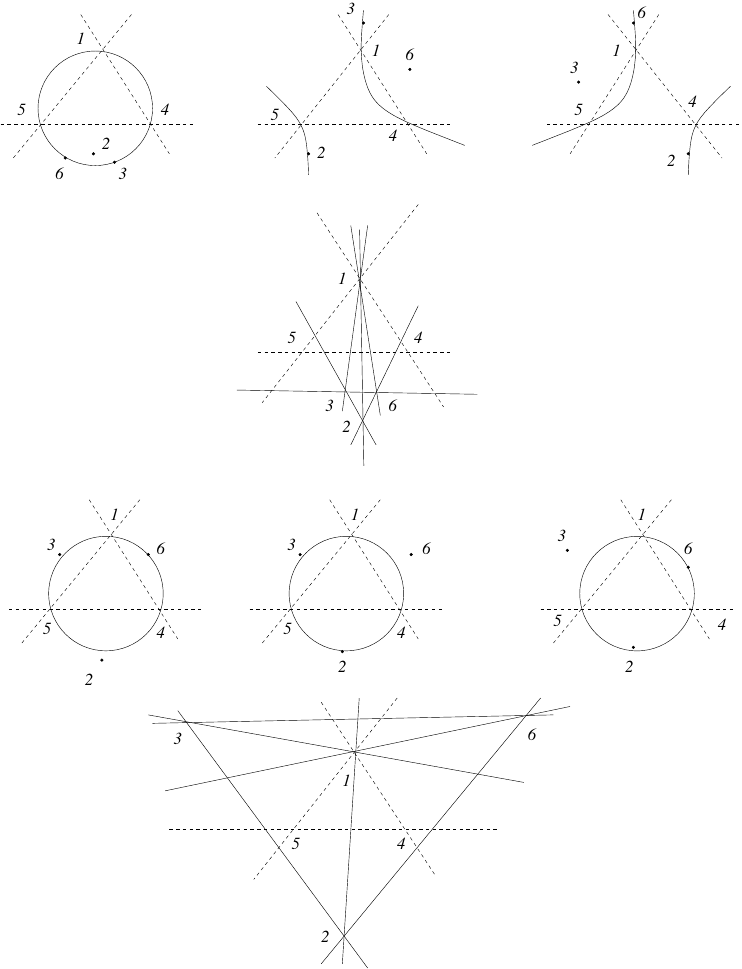}
\caption{\label{doublint} Cremona transformations for zones $A$ and $B$}  
\end{figure}

\begin{figure}
\centering
\includegraphics{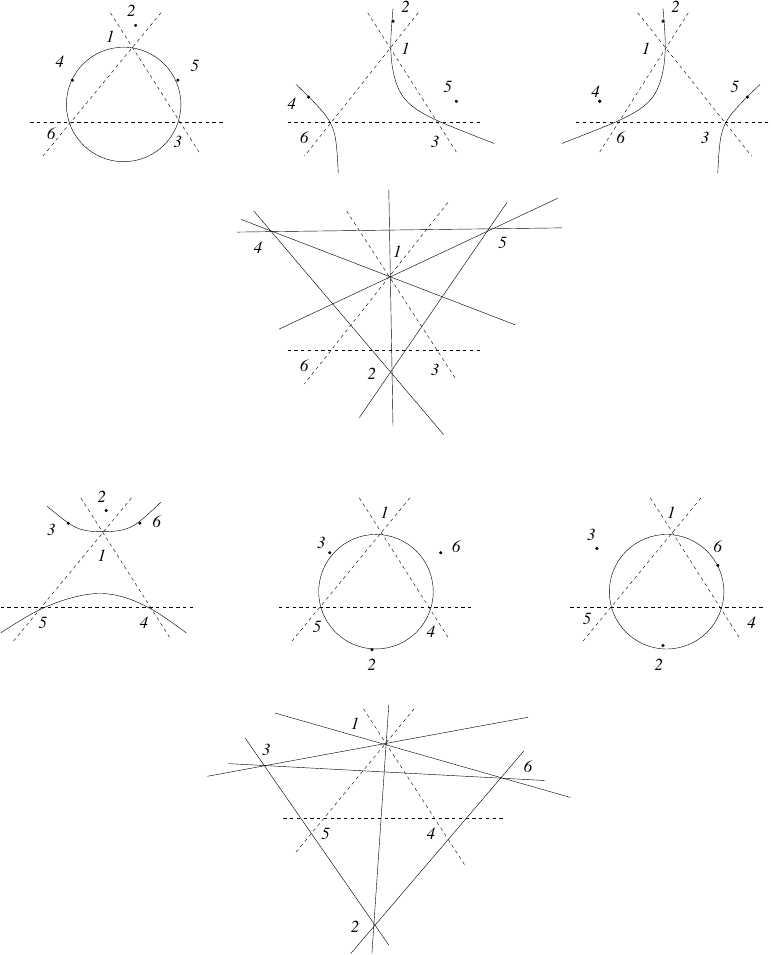}
\caption{\label{intdouble} Cremona transformations for zones $C$ and $D$}  
\end{figure}

\begin{figure}
\centering
\includegraphics{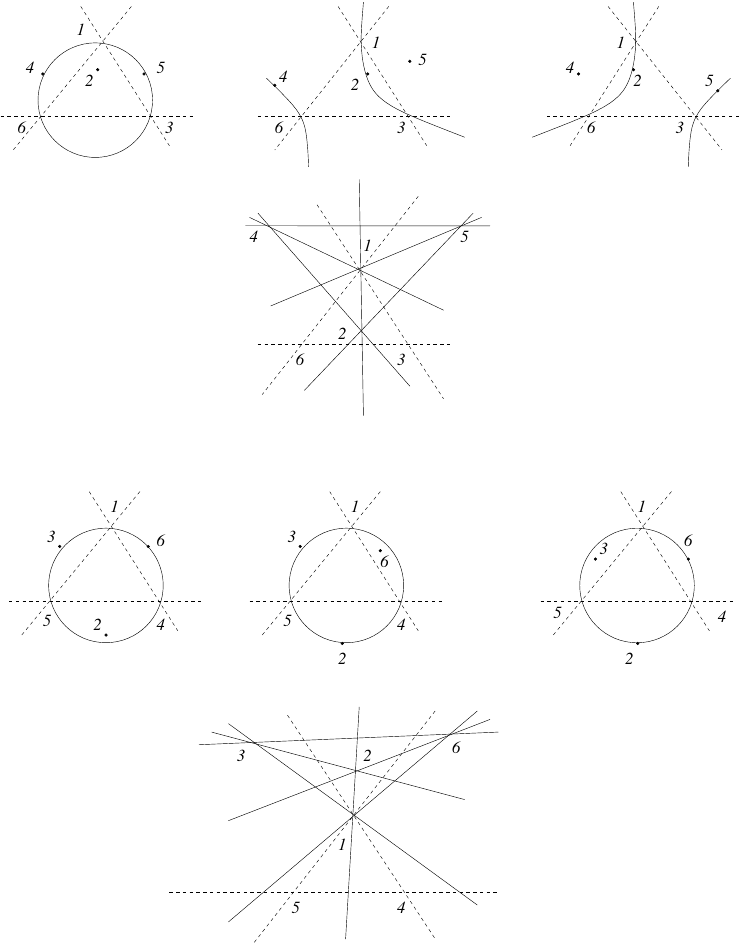}
\caption{\label{extdouble} Cremona transformations for zones $E$ and $F$}  
\end{figure}

\begin{figure}
\centering
\includegraphics{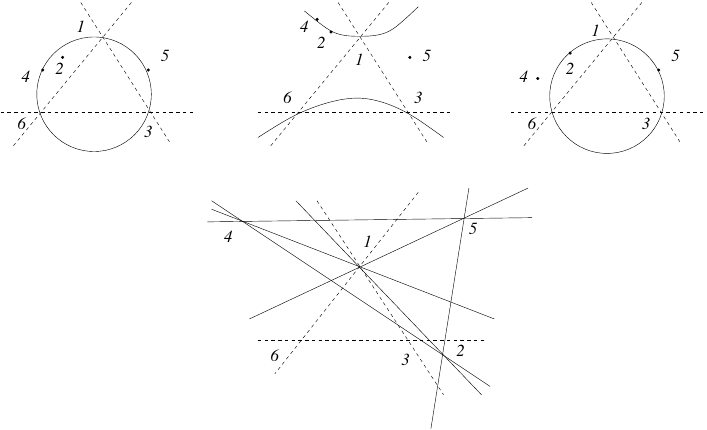}
\caption{\label{doublext} Cremona transformation for zone $G$}  
\end{figure}

\begin{figure}
\centering
\includegraphics{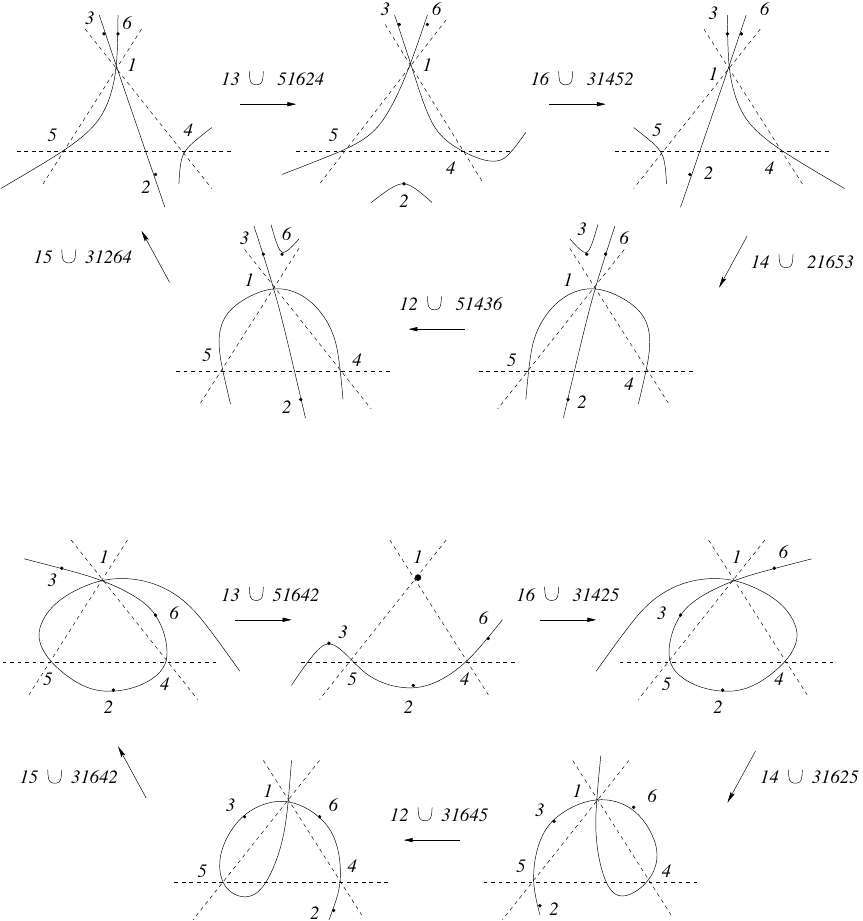}
\caption{\label{cubnod} Pencils with node in $1$ for zones $A$ and 
$B$}     
\end{figure}

\begin{figure}
\centering
\includegraphics{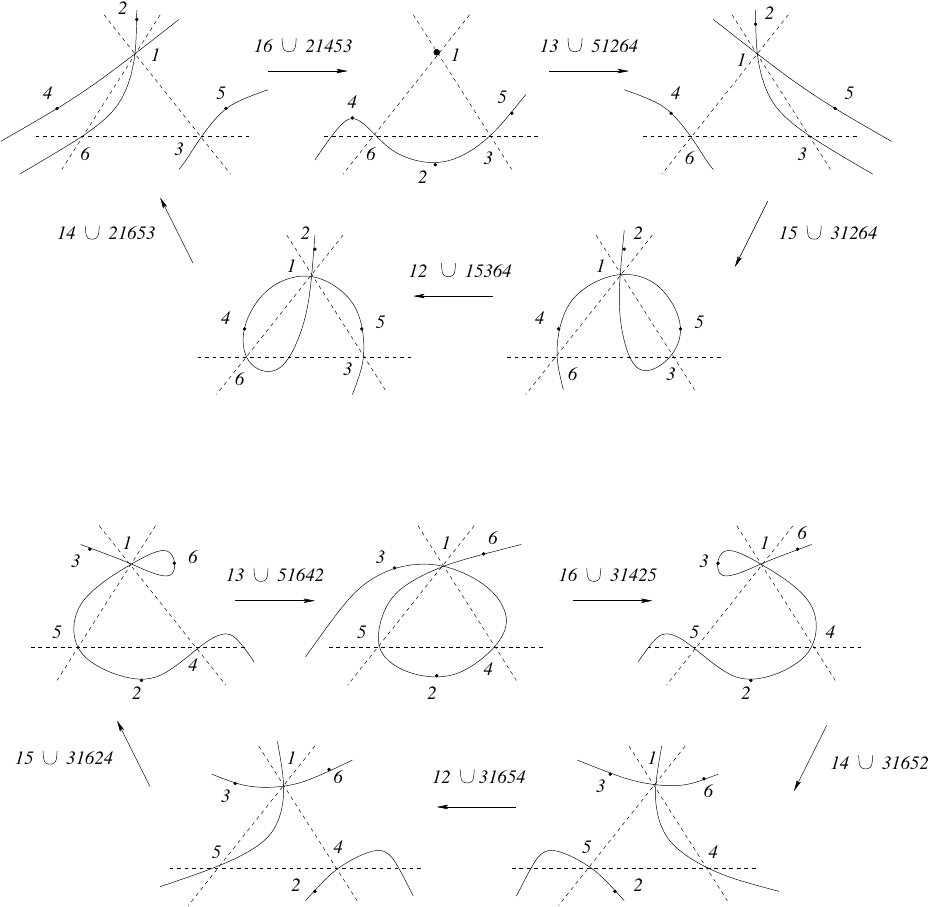}
\caption{\label{nodcub}  Pencils with node in $1$ for zones $C$ and $D$}  
\end{figure}

\begin{figure}
\centering
\includegraphics{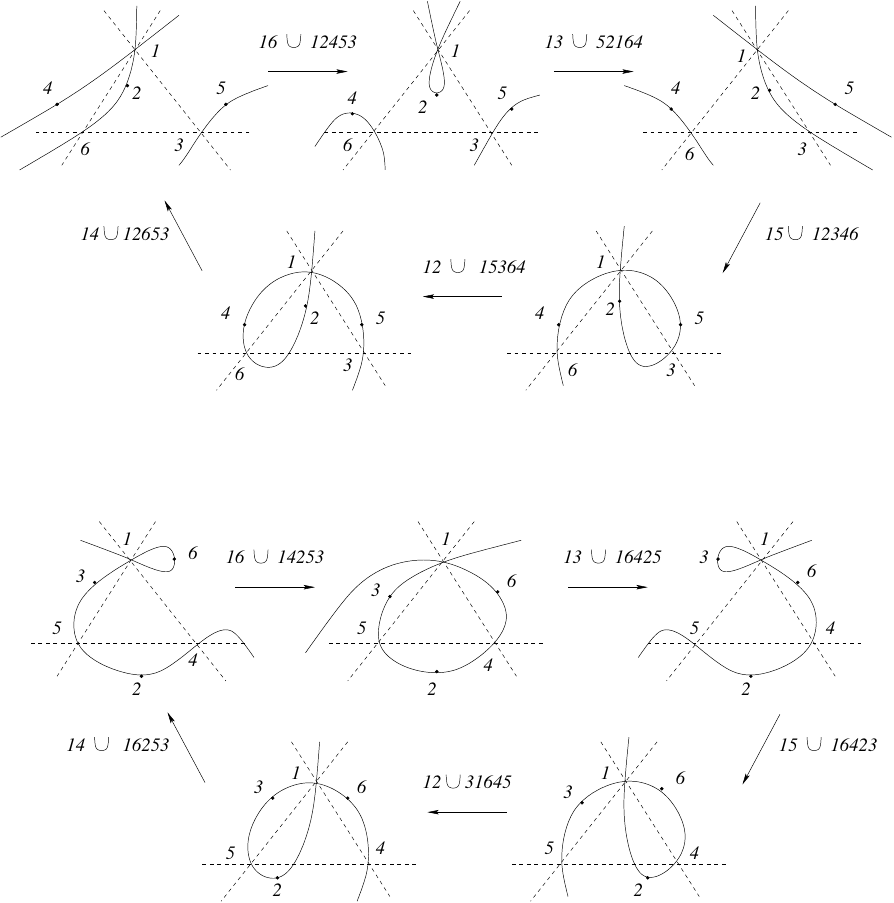}
\caption{\label{extnod} Pencils with node in $1$ for zones $E$ and 
$F$}    
\end{figure}

\begin{figure}
\centering
\includegraphics{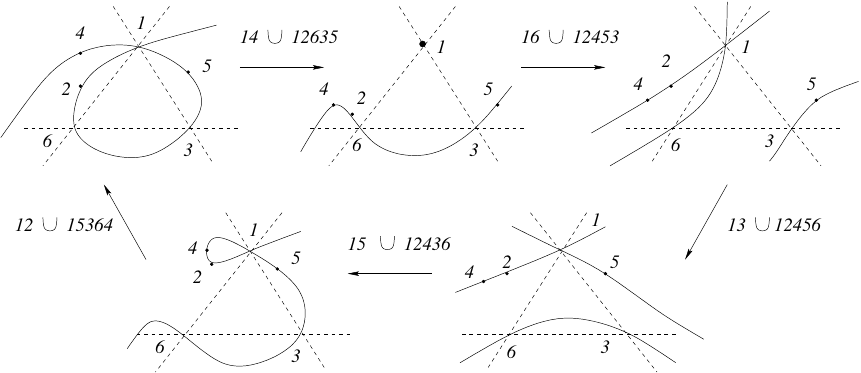}
\caption{\label{nodext} Pencil with node in $1$ for zone $G$}    
\end{figure}

Let $1 \in A$. We perform a cremona transformation $cr: (x_0; x_1; x_2) \to (x_1x_2; x_0x_2; x_0x_1)$ with base points $1, 4, 5$. Let us denote the respective images of the lines $14$, $15$, $45$ by $5$, $4$, $1$. For the other points, we keep the same notation as before $cr$. The three conics $14365$, $31452$, $61542$ are mapped onto three lines $36$, $23$ and $26$, see upper part of Figure~\ref{doublint}. After $cr$, let us consider the pencil of conics $1236$, and five particular conics of this pencil: the three double lines and the two conics passing respectively through $4$ and $5$. The cyclic ordering of these conics in the pencil is easily determined. Perform the cremona transformation back. The pencil of conics $1236$ is mapped onto the pencil of nodal cubics based at $1, \dots 6$ with node at $1$. The five particular conics are mapped  onto the five reducible cubics. 
We repeat this procedure for the various positions of $1 \in B, \dots G$, with a cremona transformation based either in $145$ or in $136$, see Figures~\ref{doublint}-\ref{doublext}. The sequences of five particular conics are displayed in Figure~\ref{conics}. The sequences of five reducible cubics for each case $A, \dots G$ are shown in Figure~\ref{seven}, the pencils are drawn in Figures~\ref{cubnod}-\ref{nodext}. More precisely, we have represented the successive types of combinatorial cubics in the five portions bounded by the reducible cubics. The upper middle pictures of each pencil $B, C, G$ are not quite correct: in the actual pencils, the portion starts and finishes with cubics having a loop containing no base point other than $1$. Inbetween, there is a pair of cuspidal cubics, and in the middle, cubics with an isolated node. We have represented only the latter. 

\begin{figure}
\centering
\includegraphics{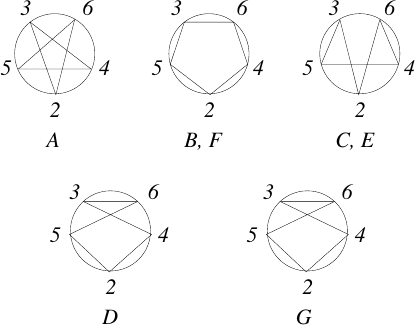}
\caption{\label{markdiag} Marked diagrams for $1 \in A, \dots G$, the circle represents the conic}    
\end{figure}

Consider any one of these seven pencils. We observe the following properties: the cyclic ordering of the points $2, \dots 6$ is the same on each of the five types of non-reducible cubic, it is given by the pencil of lines based at $1$, this follows from Bezout's theorem. The cyclic ordering of $2, \dots 6$ given by the lines of the successive reducible cubics is the same as the cyclic ordering of these points on the conic that they determine. The mutual cyclic orderings of $2, \dots 6$ given on one hand by the conic $25364$, on the other hand by the pencil of lines based at $1$, are displayed in the {\em marked diagrams\/} of Figure~\ref{markdiag}, where the circles represent the conics. 
Let us now consider together all six pencils of cubics with nodes at $1, \dots 6$.
If we take $1$ in one of the zones $E, F, G$ exterior to the conic $25364$, then we note that $2$ lies inside of the conic determined by the other five points. So, up to the action of $\mathcal{S}_6$, we may assume that $1$ is interior to a conic $25364$. There are four lists of six pencils, giving rise to four configurations, see Figure~\ref{four}. Hence, six generic points may realize four different unordered configurations $\alpha$, $\beta$, $\gamma$, $\delta$. 
For a given list, each choice of a point $1, \dots 6$ gives rise to a marked diagram.
Remove the markings (the names of the points), it turns out that all six unmarked diagrams of a list are identical, and each list corresponds to a different unmarked diagram. Note that the unmarked diagrams remain also unchanged if the roles of the circle and of the polygonal line are swapped. 
$\Box$

\begin{figure}
\begin{tabular}{c c c c c c}
$1 \in A, cr(145)$ & $12 \cup 36$ & $13462$ & $13 \cup 26$ & $16 \cup 23$ & $16532$\\
$1 \in B, cr(145)$ & $12 \cup 36$ & $31426$ & $13 \cup 26$ & $16 \cup 32$ & $25163$\\
$1 \in C, cr(136)$ & $12 \cup 45$ & $15 \cup 42$ & $41562$ & $41523$ & $14 \cup 25$\\
$1 \in D, cr(145)$ & $12 \cup 36$ & $31246$ & $13 \cup 26$ & $16 \cup 32$ & $16352$\\
$1 \in E, cr(136)$ & $12 \cup 45$ & $15 \cup 24$ & $41562$ & $41523$ & $14 \cup 25$\\
$1 \in F, cr(145)$ & $12 \cup 36$ & $34126$ & $13 \cup 26$ & $16 \cup 23$ & $32156$\\
$1 \in G, cr(136)$ & $12 \cup 45$ & $15 \cup 24$ & $41562$ & $41532$ & $14 \cup 25$
\end{tabular}
\caption{\label{conics} The seven pencils of conics}  
\end{figure}    

\begin{figure}
\begin{tabular}{c c c c c c}
$1 \in A$ & 
$12 \cup 51436$ & $15 \cup 31264$ & $13 \cup 51624$ & $16 \cup 31452$ & 
$14 \cup 21653$\\
$1 \in B$ & $12 \cup 31645$ & $15 \cup 31642$ & $13 \cup 51642$ & $16 \cup 31425$ & 
$14 \cup 31625$\\
$1 \in C$ & $12 \cup 15364$ & $15 \cup 31264$ & $13 \cup 51264$ & $16 \cup 21453$ &
$14 \cup 21653$\\
$1 \in D$ & $12 \cup 31654$ & $15 \cup 31624$ & $13 \cup 51642$ & $16 \cup 31425$ &
$14 \cup 31652$\\
$1 \in E$ & $12 \cup 15364$ & $15 \cup 12346$ & $13 \cup 12546$ & $16 \cup 12453$ & $14 \cup 12653$\\
$1 \in F$ & $12 \cup 16453$ & $15 \cup 31642$ & $13 \cup 16425$ & $16 \cup 13524$ & $14
\cup 13526$\\
$1 \in G$ & $12 \cup 14635$ & $15 \cup 12436$ & $13 \cup 12456$ & $16 \cup 12453$ & $14
\cup 12635$
\end{tabular}
\caption{\label{seven} The seven pencils of cubics}  
\end{figure}    

\begin{figure}
\begin{tabular}{c c c c c c}
$\alpha : 1 \in A$ & & & & & \\
$1 < 25364$ & $12 \cup 51436$ & $15 \cup 31264$ & $13 \cup 51624$ & $16 \cup 31452$ & 
$14 \cup 21653$\\
$2 < 51436$ & $25 \cup 31264$ & $21 \cup 25364$ & $24 \cup 21653$ & $23 \cup 51624$ &
$26 \cup 31452$\\
$3 < 51624$ & $35 \cup 31264$ & $31 \cup 25364$ & $36 \cup 31452$ & $32 \cup 51436$ &
$34 \cup 21653$\\
$4 < 21653$ & $42 \cup 51436$ & $41 \cup 25364$ & $46 \cup 31452$ & $45 \cup 31264$ &
$43 \cup 51624$\\
$5 < 31264$ & $53 \cup 51624$ & $51 \cup 25364$ & $52 \cup 51436$ & $56 \cup 31452$ &
$54 \cup 21653$\\
$6 < 31452$ & $63 \cup 51624$ & $61 \cup 25364$ & $64 \cup 21653$ & $65 \cup 31264$ &
$62 \cup 51436$\\
 & & & & & \\
$\beta : 1 \in B$ & & & & & \\
$1 < 25364$ & $12 \cup 31645$ & $15 \cup 31642$ & $13 \cup 51642$ & $16 \cup 31425$ & 
$14 \cup 31625$\\
$2 > 31645$ & $23 \cup 51642$ & $21 \cup 25364$ & $26 \cup 31425$ & $24 \cup 31625$ &
$25 \cup 31642$\\
$3 > 51642$ & $35 \cup 31642$ & $31 \cup 25364$ & $36 \cup 31425$ & $34 \cup 31625$ &
$32 \cup 31645$\\
$4 < 31625$ & $43 \cup 51642$ & $41 \cup 25364$ & $46 \cup 31425$ & $42 \cup 31645$ &
$45 \cup 31642$\\
$5 < 31642$ & $53 \cup 51642$ & $51 \cup 25364$ & $56 \cup 31425$ & $54 \cup 31625$ &
$52 \cup 31645$\\
$6 > 31425$ & $63 \cup 51642$ & $61 \cup 25364$ & $64 \cup 31625$ & $62 \cup 31645$ &
$65 \cup 31642$\\
 & & & & & \\
$\gamma : 1 \in C$ & & & & & \\
$1 < 25364$ & $12 \cup 15364$ & $15 \cup 31264$ & $13 \cup 51264$ & $16 \cup 21453$ &
$14 \cup 21653$\\
$2 > 15364$ & $21 \cup 25364$ & $25 \cup 31264$ & $23 \cup 51264$ & $26 \cup 21453$ &
$24 \cup 21653$\\
$3 > 51264$ & $35 \cup 31264$ & $31 \cup 25364$ & $32 \cup 15364$ & $36 \cup 21453$ &
$34 \cup 21653$\\
$4 < 21653$ & $42 \cup 15364$ & $41 \cup 25364$ & $46 \cup 21453$ & $45 \cup 31264$ &
$43 \cup 51264$\\
$5 < 31264$ & $53 \cup 51264$ & $51 \cup 25364$ & $52 \cup 15364$ & $56 \cup 21453$ &
$54 \cup 21653$\\
$6 > 21453$ & $62 \cup 15364$ & $61 \cup 25364$ & $64 \cup 21653$ & $65 \cup 31264$ &
$63 \cup 51264$\\
 & & & & & \\
$\delta : 1 \in D$ & & & & & \\
$1 < 25364$ & $12 \cup 31654$ & $15 \cup 31624$ & $13 \cup 51642$ & $16 \cup 31425$ &
$14 \cup 31652$\\
$2 < 31654$ & $23 \cup 51642$ & $21 \cup 25364$ & $26 \cup 31425$ & $25 \cup 31624$ &
$24 \cup 31652$\\
$3 > 51642$ & $35 \cup 31624$ & $31 \cup 25364$ & $36 \cup 31425$ & $34 \cup 31652$ &
$32 \cup 31654$\\
$4 > 31652$ & $43 \cup 51642$ & $41 \cup 25364$ & $46 \cup 31425$ & $45 \cup 31624$ &
$42 \cup 31654$\\
$5 > 31624$ & $53 \cup 51642$ & $51 \cup 25364$ & $56 \cup 31425$ & $52 \cup 31654$ &
$54 \cup 31652$\\
$6 > 31425$ & $63 \cup 51642$ & $61 \cup 25364$ & $64 \cup 31652$ & $62 \cup 31654$ &
$65 \cup 31624$
\end{tabular}
\caption{\label{four} The four lists of six pencils}  
\end{figure}

\subsection{Diagrams}
\begin{figure}
\centering
\includegraphics{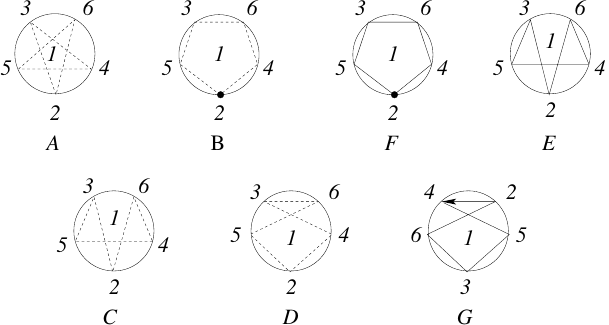}
\caption{Refined $1$-diagrams\label{refined} }    
\end{figure}

\begin{figure}
\centering
\includegraphics{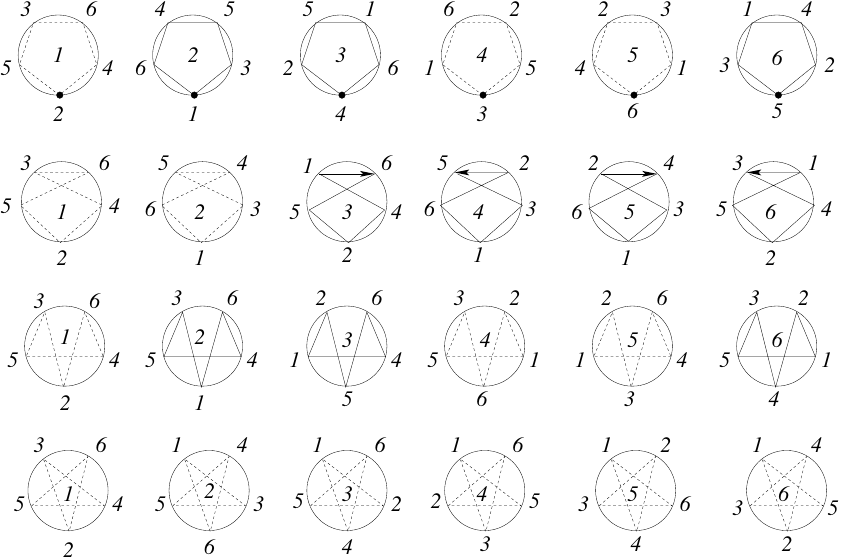}
\caption{Four configurations $\beta$, $\delta$, $\gamma$, $\alpha$ encoded each by six equivalent refined diagrams \label{equiv} }    
\end{figure}

\begin{figure}
\centering
\includegraphics{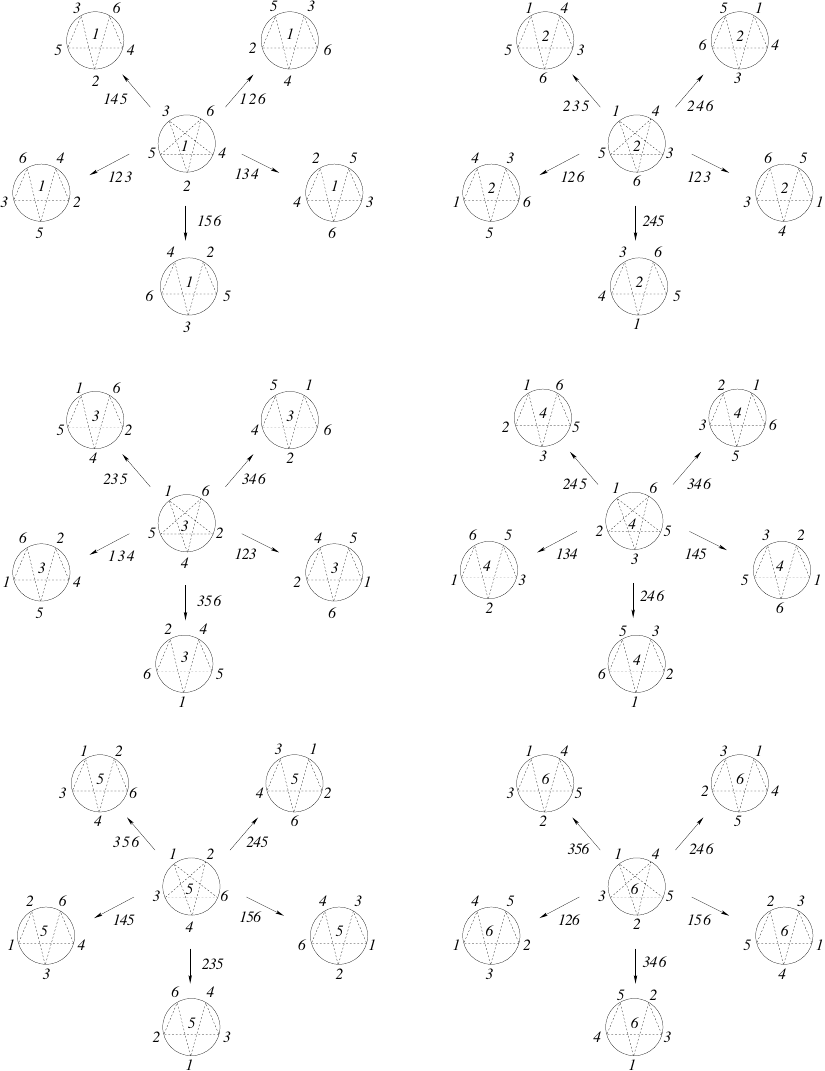}
\caption{ \label{alphamove} Crossing of walls starting from an $\alpha$-configuration, using diagrams}    
\end{figure}

\begin{figure}
\centering
\includegraphics{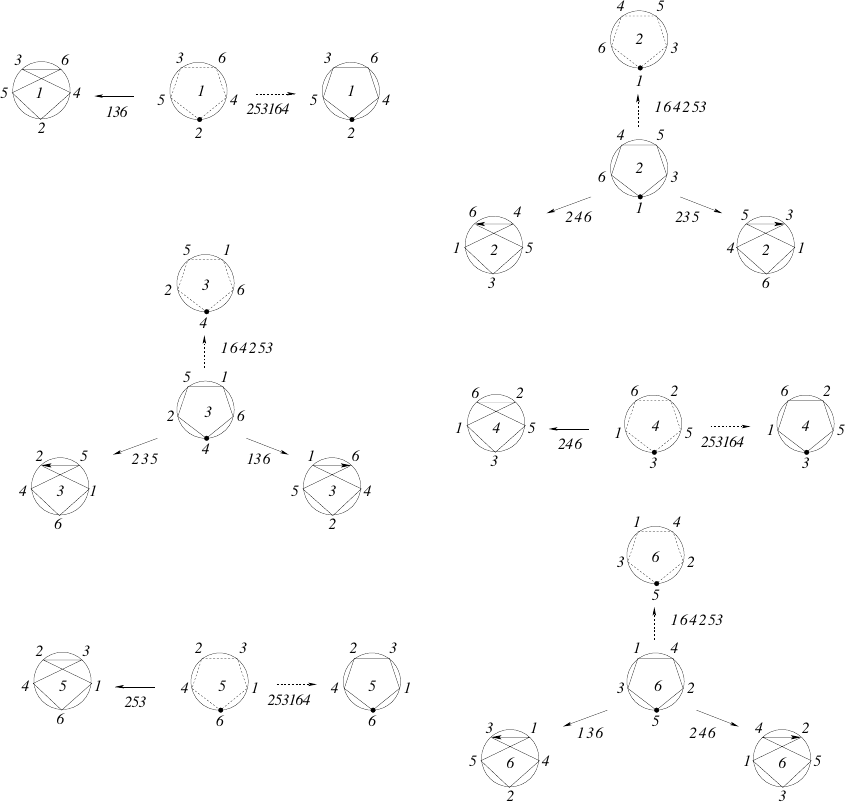}
\caption{ \label{betamove} Crossings of walls starting from a $\beta$-configuration, using diagrams}    
\end{figure}

\begin{figure}
\centering
\includegraphics{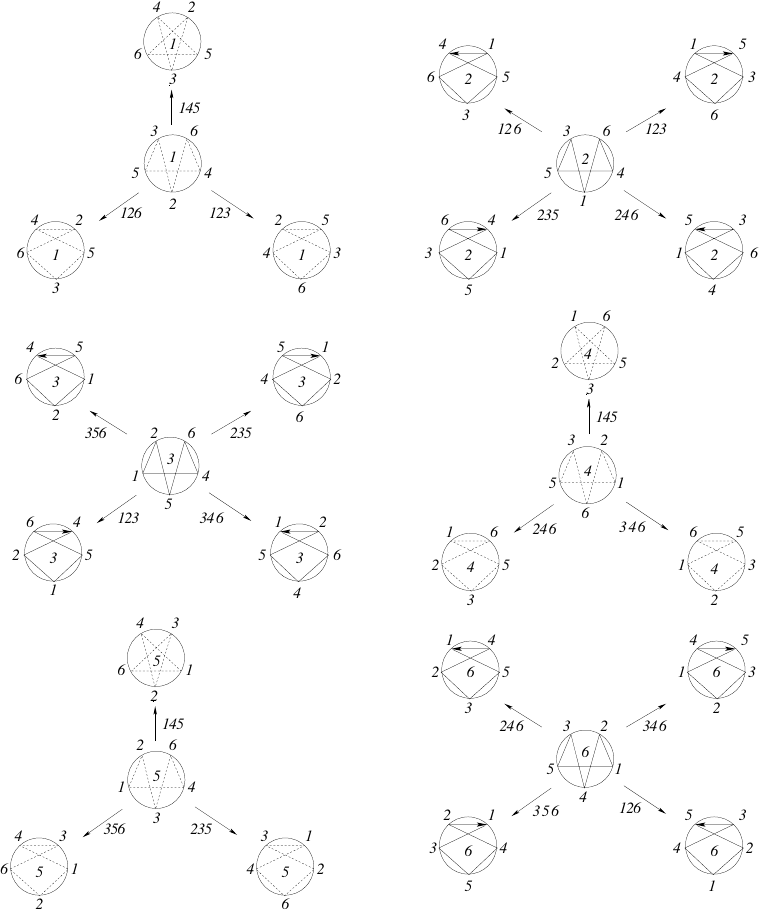}
\caption{ \label{gammamove} Crossings of walls starting from a $\gamma$-configuration, using diagrams}    
\end{figure}

\begin{figure}
\centering
\includegraphics{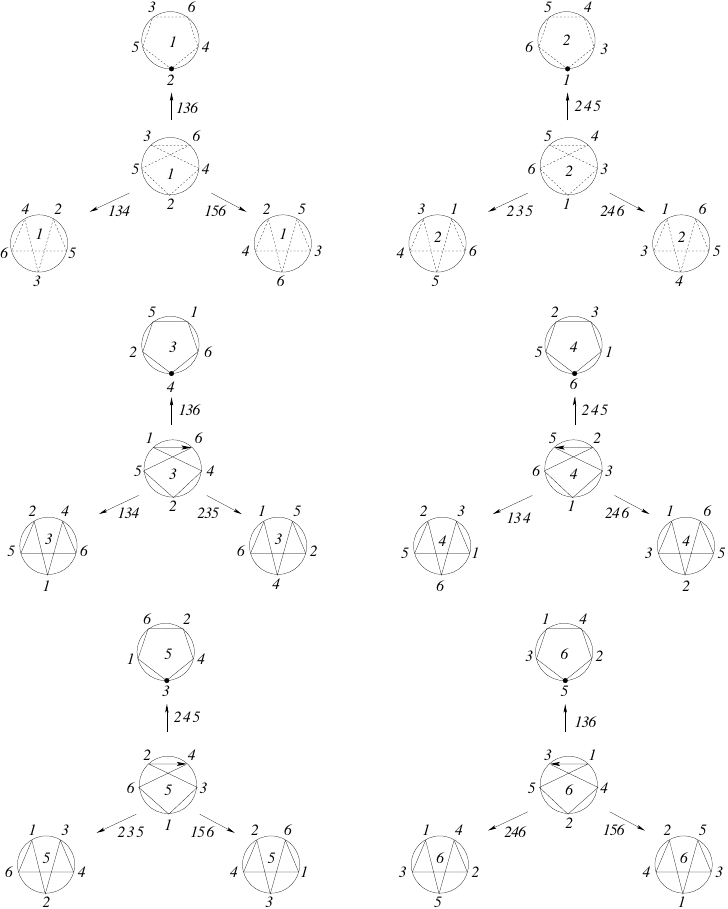}
\caption{ \label{deltamove} Crossings of walls starting from a $\delta$-configuration, using diagrams}    
\end{figure}

In the sequel, a {\em configuration\/} of points is ordered, unless otherwise explicited. The seven configurations from Figure~\ref{conf} corresponding to the choices $A, \dots G$ of the zone containing $1$ may be encoded by refining the marked diagrams from Figure~\ref{markdiag} as shown in Figure~\ref{refined}. 
First, we indicate whether $1$ lies inside or ouside of the conic $25364$
using a dotted polygonal line if $1$ is inside, and a plain polygonal line if $1$ is outside. Once this is done, the marked diagrams of $B, F$ and $G$ still correspond to several zones. The diagram of $B$ fits to the other four $B$-like zones, each of these five zones may be characterized by the point of the conic situated opposite to it, for $B$ this point is $2$. Le us add to the diagram of $B$ a dot at the point $2$. 
The same argument applies to $F$. Let us now consider the zone $G$, it is a triangle having $2$ as vertex and whose sides are supported by the lines $24$, $26$, $35$. There are in total $10$ $G$-like zones and the diagram of $G$ fits also to the one having $4$ as vertex and whose sides are supported by the lines $24$, $54$, $36$. In order to differentiate these two zones, we provide the edge $24$ in the diagram of $G$ with an arrow from $2$ to $4$. The extra point $1$ is indicated inside of 
the diagrams in Figure~\ref{refined}. Let us call $n$-diagram a diagram having $n$ as extra point.  Each configuration may actually be encoded by six diagrams, making $n = 1, \dots 6$. See  Figure~\ref{equiv}, where all six $n$-diagrams in a row are equivalent. For each of the four configurations, let one point $n$ move until it crosses a wall: a line through two other points or a conic through five other points
Figures~\ref{alphamove}-\ref{deltamove} show the corresponding changes of the $n$-diagrams for $n = 1, \dots 6$. 
We deduce the following:

\begin{proposition}
Each $\beta$-configuration is adjacent to another $\beta$-configuration via a conic-wall, and to three $\delta$-configurations via line-walls. 
Each $\delta$-configura\-tion is adjacent to two $\beta$-configurations and to four $\gamma$-configurations via line-walls.
Each $\gamma$-configuration is adjacent to six $\delta$-configurations and to one $\alpha$-configuration via line-walls.
Each $\alpha$-configuration is adjacent to ten $\gamma$-configurations via line-walls.
\end{proposition}

\begin{figure}
\centering
\includegraphics{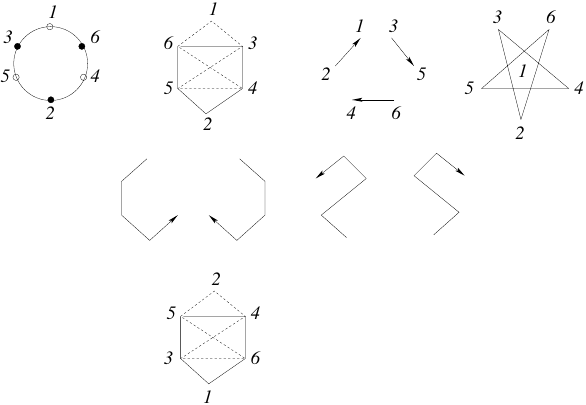}
\caption{ \label{unit} Codes for the four configurations $\beta$, $\delta$, $\gamma$, $\alpha$ from Figure 16}    
\end{figure}

It will be convenient to have a simple encoding for configurations at our disposal, so our next concern is to define a new kind of diagram, or {\em code\/}, replacing the set of six equivalent $n$-diagrams. Let us say shortly that a point $n$ is interior for the configuration if $n$ is interior to the conic determined by the other five (the $n$-diagram is dotted). The upper part of Figure~\ref{unit} shows four codes corresponding to the four configurations of Figure~\ref{equiv}. Let us explain how we defined them. In the case $\beta$, we represent the six points with six dots disposed on a circle in the natural cyclic ordering given by the convex position. The dots are colored alternatively in black and white, the white dots correspond to the interior points. The case $\delta$ offers no such evident solution, so one has to make a choice that is not entirely satisfactory. The interior points in the $\delta$-configuration of Figure~\ref{equiv} are $1$ and $2$. The two polygonal lines of their diagrams may be seen as closed paths intercepting successively five points. Let us embed these two paths in a different way, draw now one of them with a dotted and the other with a plain line. We get a graph with vertices $1, \dots 6$, having two different kinds of edges.
To recover the remaining four $n$-diagrams from this graph, note that the polygonal lines of these diagrams are oriented, with a starting point, $1$ for the $3$ and $6$-diagrams, $2$ for the $4$ and $5$-diagrams. The orderings with which the five points are met by the polygonal lines of these four diagrams are shown with arrows describing paths in the graph, as indicated in Figure~\ref{unit}. Finally, note that this encoding is not unique: the second graph drawn in the bottom of the figure represents the same configuration. Two graphs obtained from one another swapping pairs of points with a vertical symmetry represent of course also the same configuration.
For the case $\gamma$, note that the $1$-diagram is dotted and has $2$ as bottom point, whereas the $2$-diagram is plain and has $1$ as bottom point, we say that $1$ and $2$ form a pair. The other four points may be similarly distributed in two pairs. Let us encode this with a graph having the six points as vertices, and three edges endowed with arrows as shown in Figure~\ref{unit}, each arrow connects an interior point to its associated exterior point. The original $n$-diagrams may be deduced easily from this graph.
For the case $\alpha$, let us simply observe that the $n$-diagrams may be deduced easily from one another. For example, start with the $1$-diagram. In this diagram, $1$ and $2$ are separated by the branch $54$ of the star. To get the $2$-diagram, it suffices to swap $1$ with $2$ and $5$ with $4$.  To encode this configuration, we will simply use any one of the $n$-diagrams, with the circle removed.

\begin{figure}
\centering
\includegraphics{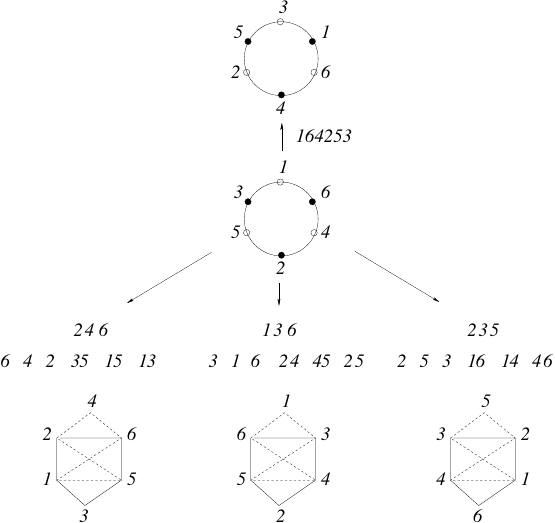}
\caption{ \label{movebeta} Crossings of walls starting from a $\beta$-configuration, using codes}    
\end{figure}

\begin{figure}
\centering
\includegraphics{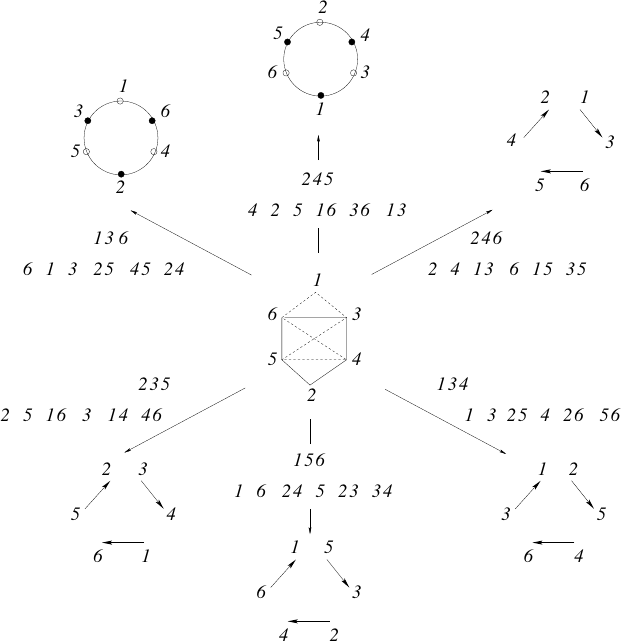}
\caption{ \label{movedelta} Crossings of walls starting from a $\delta$-configuration, using codes}    
\end{figure}

\begin{figure}
\centering
\includegraphics{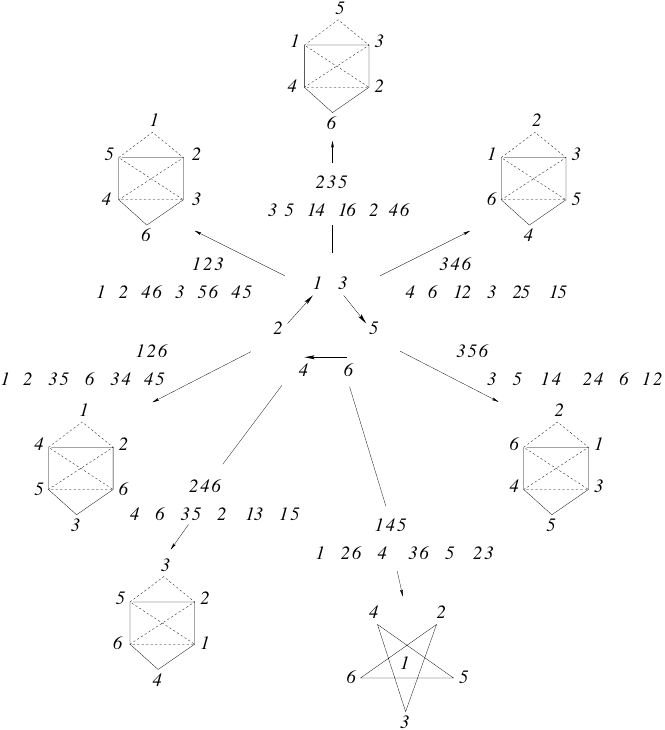}
\caption{ \label{movegamma} Crossings of walls starting from a $\gamma$-configuration, using codes}    
\end{figure}

\begin{figure}
\centering
\includegraphics{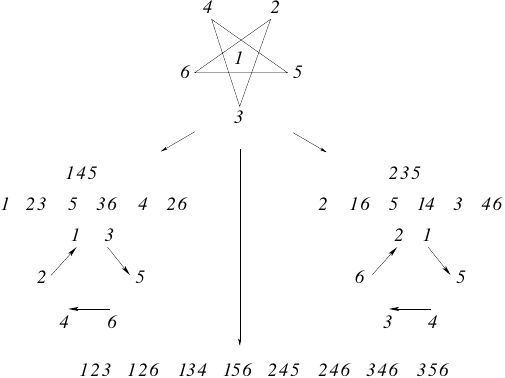}
\caption{ \label{movealpha} Crossing of walls starting from an $\alpha$-configuration, using codes}    
\end{figure}

\begin{figure}
\centering
\includegraphics{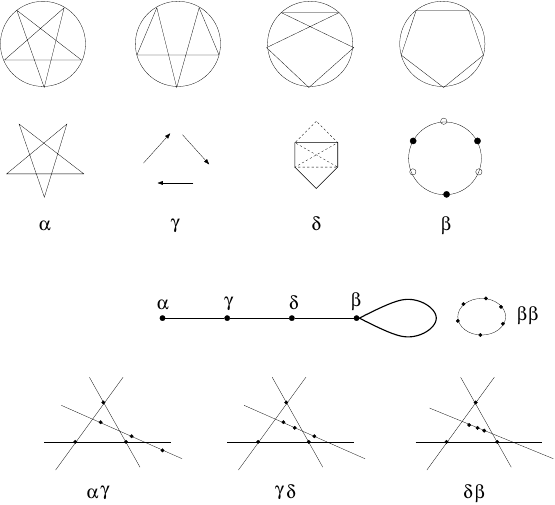}
\caption{\label{diagcodes} Adjacency graph for six unordered points, stratification by lines and conics}    
\end{figure}

Figures~\ref{movebeta}-\ref{movealpha} show all possible crossings of walls, starting from the four configurations, using now codes instead of diagrams. Near each arrow corresponding to a wall we indicate the line passing through three points or the conic passing through six points, with a triple or a cyclically ordered sextuple.
In Figure~\ref{movealpha}, the list of ten lines is written, but for place reason, we drew the codes of only two adjacent $\gamma$-configurations.
A line-wall is described by the cyclic ordering with which the line meets: the three points and the three lines determined by the other three points. 
The quotient space $(\mathbb{R}P^2)^6/\mathcal{S}_6$ has an algebraic variety structure as orbit space of a finite group action. Let us endow it with the natural stratification given by the alignment of three points, it has four cameras and three line-walls. Refine the stratification with the conics, one conic-wall appears inside of one camera. The adjacency graph, obtained assigning a vertex to each camera and an edge to each wall, was first described in \cite{sf1}-\cite{sf2}. (In these papers, Finashin used a different approach, considering actually arrangements of lines dual to configurations of points). 
The three line-walls are: $\alpha \gamma: (12, 4, 13, 5, 23, 6)$, $\gamma \delta: (12, 4, 5, 13, 23, 6)$ and $\delta \beta: (12, 13, 23, 4, 5, 6)$ (where $4, 5, 6$ stand for the three aligned points), the conic-wall is $\beta \beta: 123456$, see
Figure~\ref{diagcodes}. 
Each of the four configurations $\beta$, $\delta$, $\gamma$, $\alpha$ from Figure~\ref{unit} is preserved by some subgroup $G$  of $\mathcal{S}_6$, which may be found directly from the code. One has: 

$G(\beta) = \{ id, (3 6)(4 5), (1 4)(2 3), (1 5)(2 6), (1 4 5)(2 3 6), (1 5 4)(2 6 3) \} = D_3$,

$G(\delta) = \{ id, (3 6)(4 5), (1 2)(3 4 6 5), (1 2)(3 5 6 4) \} = \mathbb{Z}/4$,
 
$G(\gamma) = \{ id, (3 6)(4 5)$, $(3 2)(1 5), (1 4)(2 6)$, $(1 5 4)(2 3 6), (2 6 3)(1 4 5) \} = D_3$. 

The group $G(\alpha)$ is a icosahedral group (isomorphic to $A_5$), it has $60$ elements: $id$, $24$ elements of order two (conjugacy class of $(1 2)(4 5)$), $20$ of order five (conjugacy class of $(2 5 3 6 4)$), and $15$ of order three (conjugacy class of $(1 4 2)(3 5 6)$). 
The action of $G(\alpha)$ on the points is transitive. For each of the other three configurations, the points are distributed in two orbits, so there are in total seven pairs (configuration, orbit): $a$ (for $\alpha$); $b1$:$\{ 2, 3, 6 \}$, $b2$:$\{ 1, 4, 5 \}$ ($\beta$); $c1$:$\{ 2, 3, 6 \}$, $c2$:$\{ 1, 5, 4 \}$ ($\gamma$), $d1$:$\{ 1, 2 \}$, $d2$:$\{ 3, 4, 5, 6 \}$ ($\delta$).
In the quotient space, the walls adjacent to a configuration correspond to the orbits of the adjacent triples and eventual sextuple under the action of its monodromy group.
For example $G(\delta)$ gives rise to two orbits of triples: $\{ 136, 245 \}$ ($\delta \beta$) and $\{ 246, 134, 156, 235 \}$ ($\delta \gamma$), see Figure~\ref{movedelta},


A configuration describes the mutual position of points with respect to lines and conics, we could call it a $Q$-configuration. If we leave aside the information about the conics, we speak of an $L$-configuration. A (generic) $L$-configuration may be described in the same way as a $Q$-configuration, except that one considers only the $31$ zones determined by the lines instead of the $36$ zones determined by the lines and the conic.
A $L$ (resp. $Q$)-{\em deformation\/} is a generic deformation in the space $(\mathbb{R}P^2)^6$ stratified with the lines only (resp. with the lines and the conics). 
Let us encode $\alpha$, $\delta$, $\gamma$, $\beta$ now as $L$-configurations. The codes for the first three stay unchanged. For $\beta$, one has to remove the colors black and white of the dots. 
Each $L$-configuration is preserved by some subgroup $G'$ of $\mathcal{S}_6$. One has  $G'(\beta) = D_6$, note that this group acts transitively on the points. For each of the other three configurations, one has $G' = G$. There are in total six pairs (configuration, orbit): $a$, $b$, $c1$, $c2$, $d1$, $d2$. In $(\mathbb{R}P^2)^6$ stratified with the lines only, two generic elements realizing the same topological type (configuration) are rigidly isotopic \cite{sf1}. The same holds when one refines the stratification with conics \cite{z}-\cite{fz}. So, the group $G'$ (resp. $G$) associated to a $L$ (resp. $Q$)-configuration is the {\em monodromy group\/} of $(1, \dots 6) \in (\mathbb{R}P^2)^6$, i.e. the  subgroup of $\mathcal{S}_6$ formed by the permutations of these points realized by $L$ (resp. $Q$)-deformations.  
Two sextuples of points $(P_1, \dots P_6)$ and $(P'_1, \dots P'_6)$ realizing a wall 
are rigidly isotopic if (and only if) they have the same topological type,
let us explicit this. Assume first that $P_4, P_5, P_6$ lie on a line $L$ and $P'_4, P'_5, P'_6$ lie on a line $L'$. Either set of six points determines four lines in total, let $L_{ij}$ be the line through $P_i, P_j$ and $L'_{ij}$ be the line through $P'_i, P'_j$. A projective transformation $\sigma$ maps $(L'_{12},
L'_{13}, L'_{23}, L')$ onto $(L_{12}, L_{13}, L_{23}, L)$. Note that $\sigma$ is realizable with a deformation of the whole plane, as $PGL(3, \mathbb{R})$ is connected. One has $\sigma(P'_i) = P_i$, $i = 1, 2, 3$. There exists a
further deformation, shifting $\sigma(P'_i)$, $i = 4, 5, 6$ onto $P_i$ along $L$. 
Assume now that $(P_1, \dots P_6)$ lie in this ordering on a conic $C$, and $(P'_1, \dots P'_6)$ on a conic $C'$.
A projective transformation $\tau$ maps $C'$ onto $C$.
A further deformation along $C$ shifts $(\tau(P'_1), \dots \tau(P'_6))$ onto
$(P_1, \dots P_6)$.

Let us mention finally the monodromy groups of the walls:

$G(24, 6, 1, 3, 25, 45) = \{ id, (3 6)(4 5) \}$, 

$G(34, 1, 6, 24, 5, 23) = \{ id \}$, 

$G(23, 1, 26, 4, 36)$ = $\langle (3 6)(4 5), (2 3 6)(1 5 4) \rangle$, 

$G(123456)$ = $D_6 = \langle (1 2)(4 5), (1 2 3 4 5 6) \rangle$.

The conic-wall connects two $\beta$-configurations that have the same monodromy group $G = \{ id, (2 6)(3 5), (1 5)(2 4), (1 3)(4 6), (1 3 5)(2 4 6), (1 5 3)(2 6 4) \}$
Each of the line-walls is the biggest common subgroup of the groups associated to the two adjacent configurations. (In the case of $\beta \delta$, this is true with both $G(\beta)$ and $G'(\beta)$).

\section{Configurations of seven points}
\subsection{Fourteen configurations}

A configuration of seven points $1, \dots 7$ will be encoded with the list of codes for the subconfigurations $\hat 1, \dots \hat 7$. Assume that $1, \dots 6$ are disposed in this ordering on a conic. The non-generic configuration $\hat 7$ may be encoded with a circle passing successively through six points $1, \dots 6$, it is a $\beta$-code from which the colors black and white were removed. (To avoid a confusion with the code of an $L$-configuration, one may add a letter $c$ for conic inside.)
The mutual cyclic orderings of $1, \dots, 6$ given respectively by the conic and the pencil of lines based at $7$ may be described with a {\em conic-diagram\/}: a closed polygonal line with six vertices, inscribed in a circle. It is easily seen that there are eleven admissible such unmarked conic-diagrams. 

\begin{proposition}
Seven points $1, \dots 7$ in $\mathbb{R}P^2$ with six of them 
on a conic, but otherwise generic, may realize eleven different unordered configurations.

Seven generic points $1, \dots 7$ in $\mathbb{R}P^2$ that lie in convex position may realize eleven different unordered configurations. Up to the action of $\mathcal{S}_7$, there are correspondingly eleven lists of seven nodal combinatorial cubics.
\end{proposition} 
\begin{figure}
\centering
\includegraphics{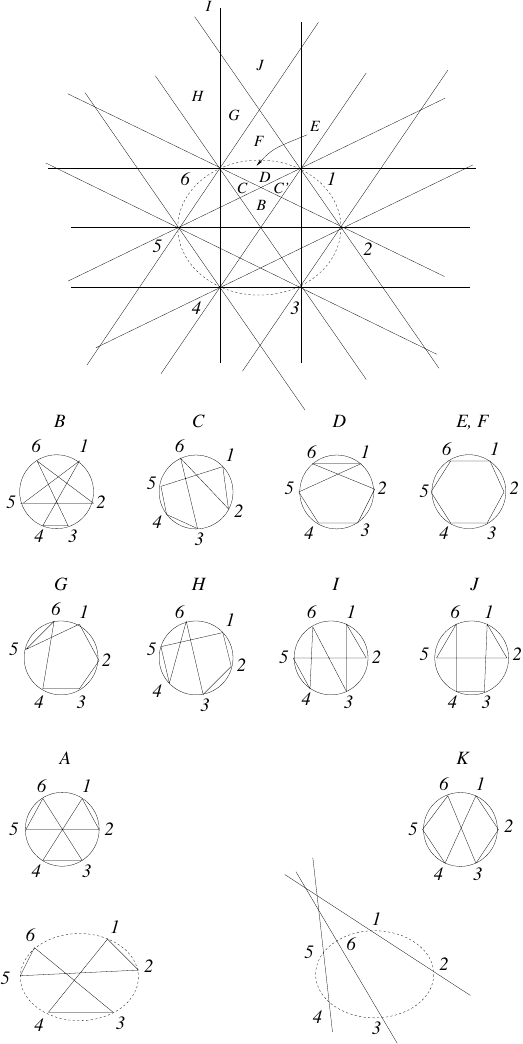}
\caption{ \label{hexagon} The eleven configurations $A, \dots K$ with six coconic points}    
\end{figure}

\begin{figure}
\centering
\includegraphics{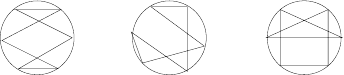}
\caption{ \label{impossible} The three non-realizable unmarked conic-diagrams}    
\end{figure}

\begin{figure}
\centering
\includegraphics{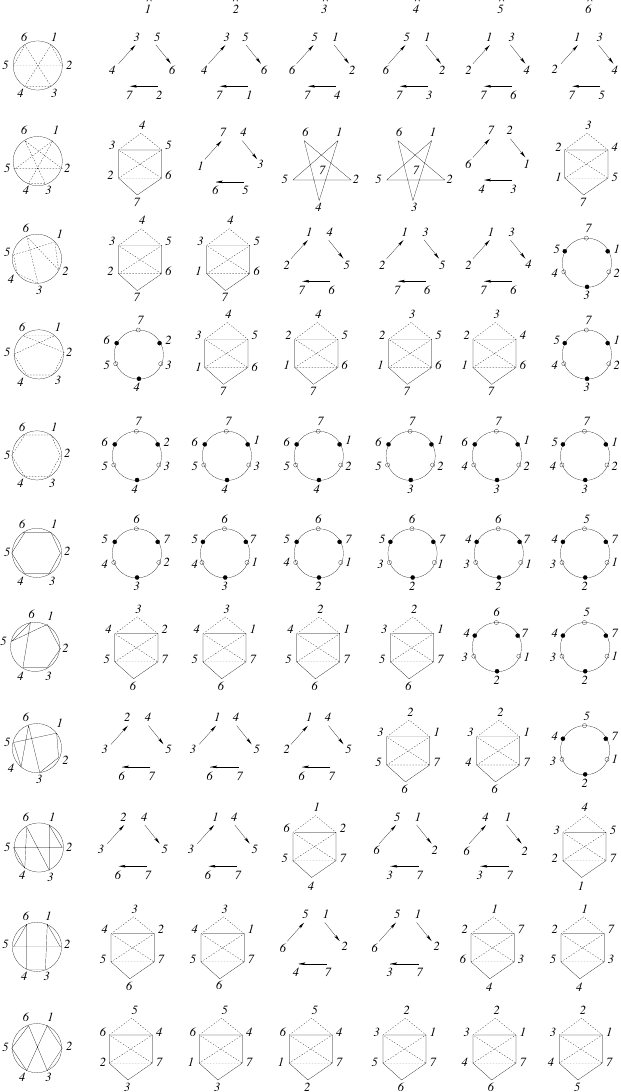}
\caption{ \label{coconicone} Conic-diagrams and codes $\hat 1, \dots \hat 6$ for
the configurations $A, \dots K$}    
\end{figure}

\begin{figure}
\centering
\includegraphics{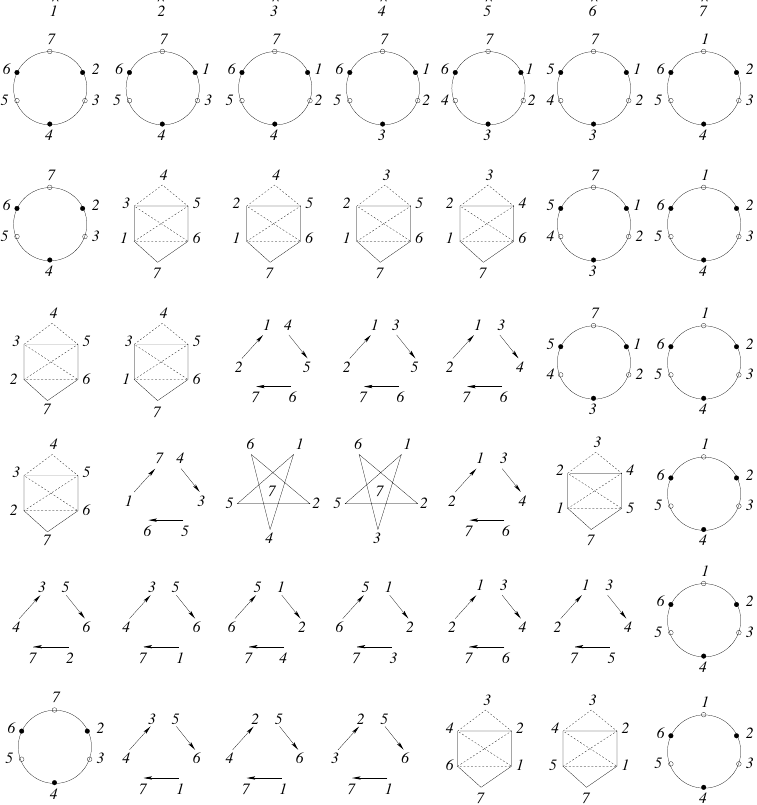}
\caption{ \label{convone} Codes of $(E, 6)$, $(D, 6)$, $(C, 6)$, $(B, 6)$, $(A, 6)$ and $(C', 6)$}    
\end{figure}

\begin{figure}
\centering
\includegraphics{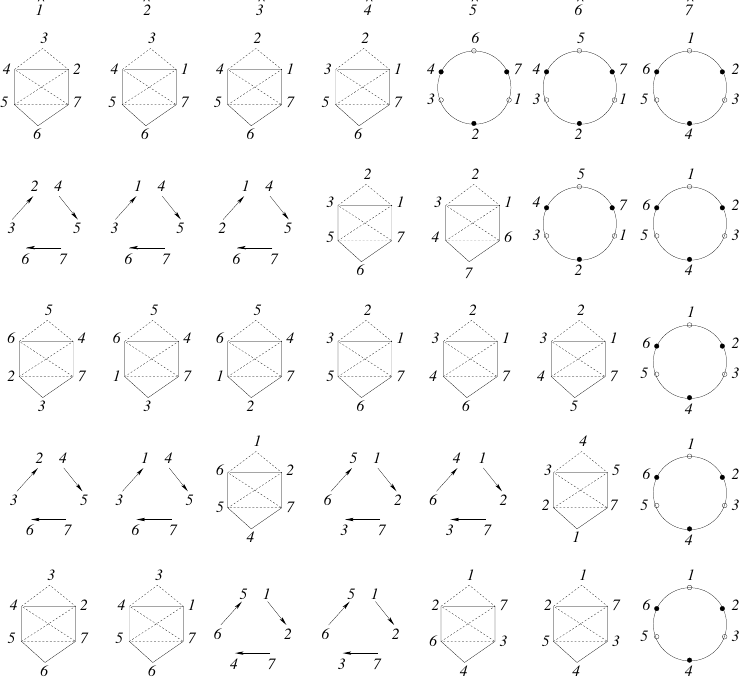}
\caption{ \label{convtwo} Codes of $(G, 6)$, $(H, 6)$, $(K, 6)$, $(I, 6)$ and $(J, 6)$}    
\end{figure}

\begin{figure}
\centering
\includegraphics{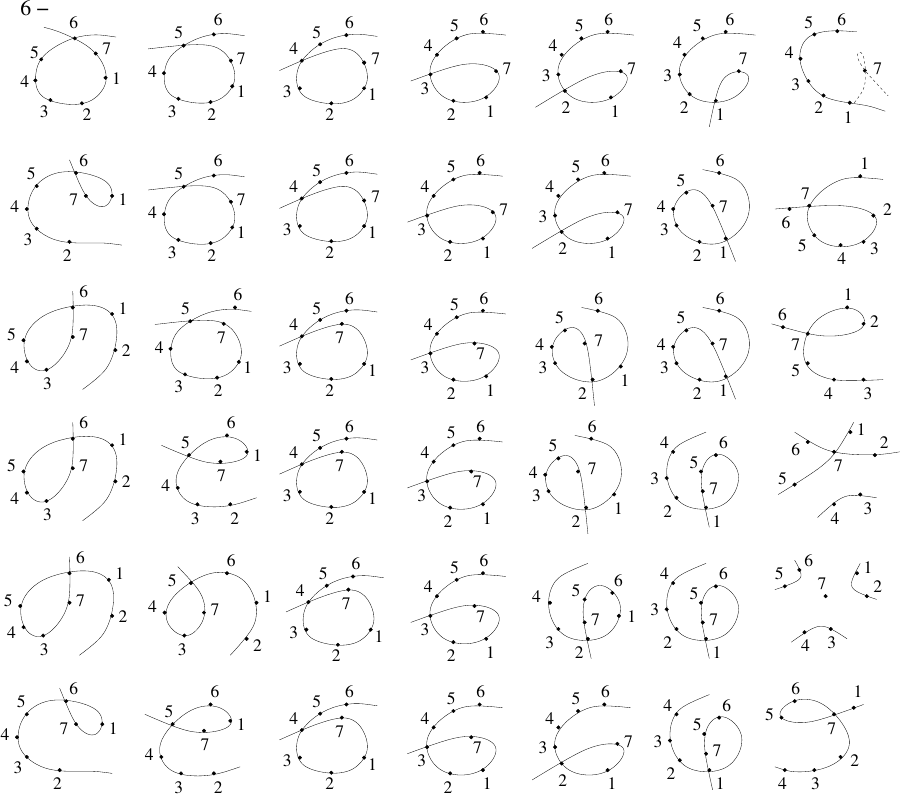}
\caption{ \label{sevenfirst} Lists of seven nodal cubics for $(E, 6)$, $(D, 6)$, $(C, 6)$, $(B, 6)$, $(A, 6)$ and $(C', 6)$}    
\end{figure}

\begin{figure}
\centering
\includegraphics{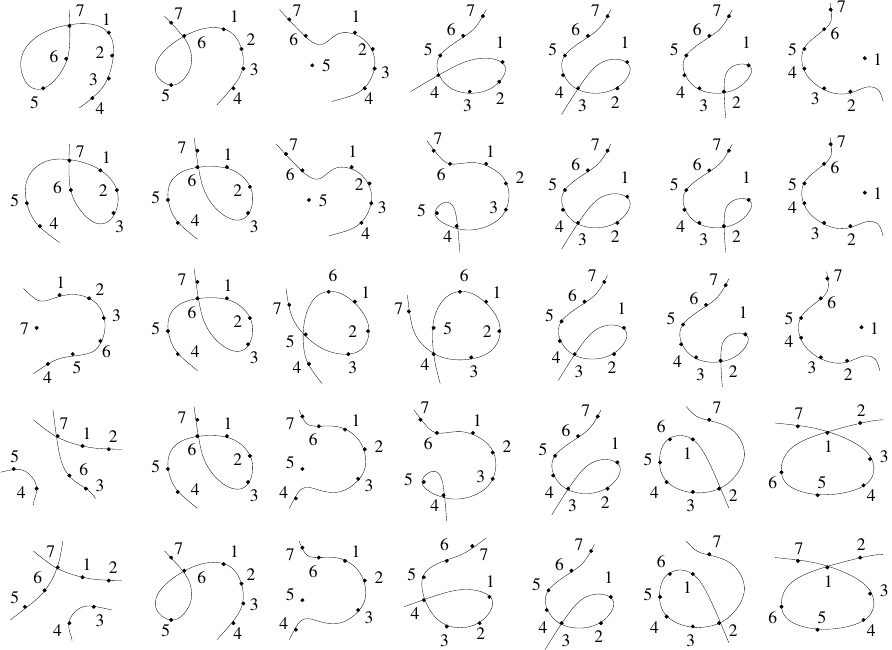}
\caption{ \label{secondseven} Lists of seven nodal cubics for $(G, 6)$, $(H, 6)$, $(K, 6)$, $(I, 6)$ and $(J, 6)$}    
\end{figure}

{\em Proof:\/} 
Let us consider for a start seven points $1, \dots 7$ such that $1, \dots 6$ lie in this ordering on a conic, and draw this conic as an ellipse in some affine plane. The point $7$ lies in a zone bounded by some of the $15$ lines determined by $1, \dots 6$. The lines $36$, $14$, $25$ give rise to six sectors, containing each one edge of the hexagonal convex hull of the six points. Up to cyclic permutation of $1, \dots 6$, we may assume that $7$ lies in the sector containing the edge $61$.  Let us move the six points keeping them coconic until three triplets of lines become concurrent, as shown in the upper part of Figure~\ref{hexagon}. Move $7$ along so that it does not cross any of the $15$ lines nor the conic. Note that the zone containing $7$ may be a triangle that shrinks in the end to a triple point. If $7$ is not in such a vanishing zone, one may assume up to the symmetry $(6 1)(5 2)(4 3)$ that $7$ lies in the end in one of the nine zones $B, \dots J$ of Figure~\ref{hexagon}. Note that the configuration realized by the seven points is preserved all along the motion. If $7$ is in a vanishing triangle, this triangle is either $A$ or $K$, see the bottom part of Figure~\ref{hexagon} (up to the action of $(6 1)(5 2)(4 3)$ for $K$). There are thus eleven unordered configurations of seven points with six of them coconic, name them $A, \dots K$ after the zone containing $7$. Note that these eleven zones give rise in total to eight unmarked conic-diagrams, the other three that turn out to be unrealizable are shown in Figure~\ref{impossible}. 
Figure~\ref{coconicone} shows for each zone $A, \dots K$: the unmarked conic-diagram (refined with either a dotted or a plain polygonal line depending on whether $7$ lies inside or outside of the conic), and the six codes $\hat 1, \dots \hat 6$. Let now $1, \dots 7$ be seven generic points, such that six of them lie in convex position. We may assume that the seventh point (extra point) is either inside or outside of all six conics determined by the first six. Let indeed $1, \dots 6$ lie in convex position, $6$ being outside of the conic $12345$, and let $7$ lie between two of the six conics. If the seven points lie in convex position (the position of $7$ in the cyclic ordering is arbitrary), consider one of the two conics adjacent to $7$, let $n \in 1, \dots 6$ be the point that doesn't lie on this conic, $n$ may be taken as extra point; otherwise, up to the action of $D_3$ on $1,\dots, 6$, one of the following two situations occurs: $1, 2, 3, 4, 5, 7$ lie in convex position and $6$ is outside of all six conics they determine, or $2, 3, 4, 5, 6, 7$ lie in convex position and $1$ lies inside of all six conics they determine. 
Let now $1, \dots 7$ lie in convex position, and $7$ be either inside or outside of all
six conics determined by $1, \dots 6$. One may move the seven points until the first six become coconic, preserving the configuration all along.
So, up to the action of $\mathcal{S}_7$, any configuration of seven points with six
of them in convex position may be obtained from one of the eleven non-generic configurations $A, \dots K$ by moving $6$ away from the conic $12345$.  Let $X$ be one of the zones in Figure~\ref{hexagon}. Denote by $(X, 6)$ the configuration obtained from $X$ moving $6$ to the outside of the conic $12345$, and $(X, 6')$ the configuration obtained moving $6$ to the inside. 
Let us list the pairs of equivalent configurations, along with the elements of
$\mathcal{S}_7$ mapping one onto the other:
$(X, 6)$ and $(X, 6')$ with $X \in \{ A, B, D, E, F, J \}$ are swapped by the symmetry $(6 1)(5 2)(4 3)$, this symmetry swaps also $(C', 6)$ with $(C, 6')$.  The symmetry $(4 2)(7 6)(1 5)$ swaps $(D, 6)$ with $(G, 6')$, and $(C, 6)$ with $(H, 6')$. The symmetry $(6 3)(1 4)(5 2)$ swaps $(I, 6)$ with $(I, 6')$, and $(K, 6)$ with $(K, 6')$.  
The cyclic permutation $(1 2 3 4 5 6 7)$ maps $(E, 6)$ onto $(F, 6')$.  
There are thus eleven different unordered configurations of seven points with six in convex position. Let us name for first each of them after some representant for the equivalence class. Later on, we will introduce a more canonical encoding. 
Each of the zones $X \in \{A, B, I, J, K\}$ gives rise to one unordered configuration, say $(X, 6)$, the pair $E, F$ gives rise to one, say $(E, 6)$, the pair $D, G$ gives rise to two, say $(D, 6)$ and $(G, 6)$ and the pair $C, H$ gives rise to three, say $(C, 6)$, $(C', 6)$ and $(H, 6)$. Figure~\ref{convone} and Figure~\ref{sevenfirst} display respectively the codes and the lists of seven nodal cubics with nodes at $1, \dots 7$ for the configurations $(E, 6)$, $(D, 6)$, $(C, 6)$, $(B, 6)$, $(A, 6)$, $(C', 6)$.  Figure~\ref{convtwo} and Figure~\ref{secondseven} display the codes and the lists of cubics for the configurations $(G, 6)$, $(H, 6)$, $(K, 6)$, $(I, 6)$, $(J, 6)$. Note that $(E, 6)$ is the list denoted by $6-$ in \cite{fi4}. Let us explain how we find out these lists of cubics. To get a cubic with node at a given point, say $1$, it suffices to perturb the reducible cubic $17 \cup 123456$, moving one point, say $6$, away from the conic $12345$ in the appropriate direction. Note that crossing a line $lnm$
induces the change of the four subcodes $\hat p$, $p \not= l, n, m$, and of the three cubics with respective node at $l, n, m$. $\Box$

\begin{figure}
\centering
\includegraphics{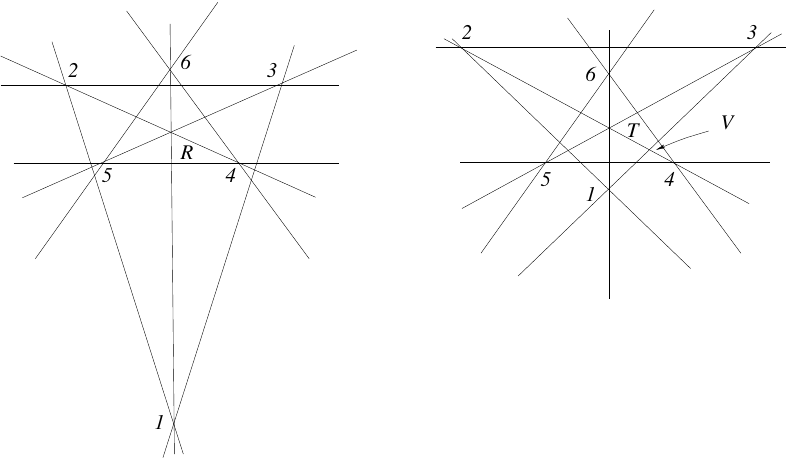}
\caption{ \label{nonconvone} Three new configurations: $7 \in R$, $7 \in T$, $7 \in V$}    
\end{figure}

\begin{figure}
\centering
\includegraphics{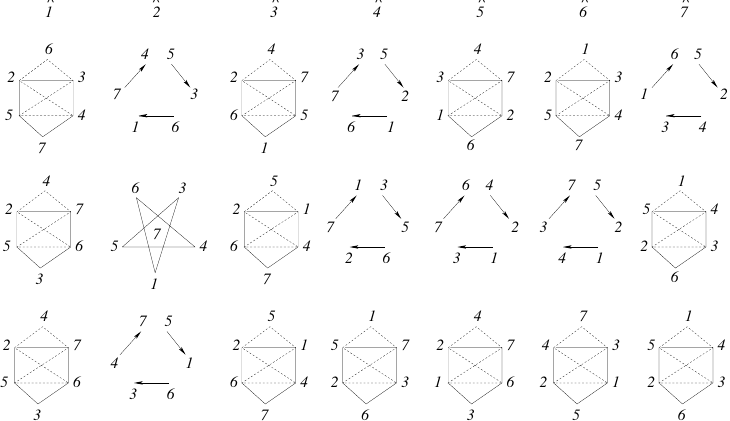}
\caption{ \label{nonconvtwo} Codes of $R$, $T$, $V$}    
\end{figure}

\begin{figure}
\centering
\includegraphics{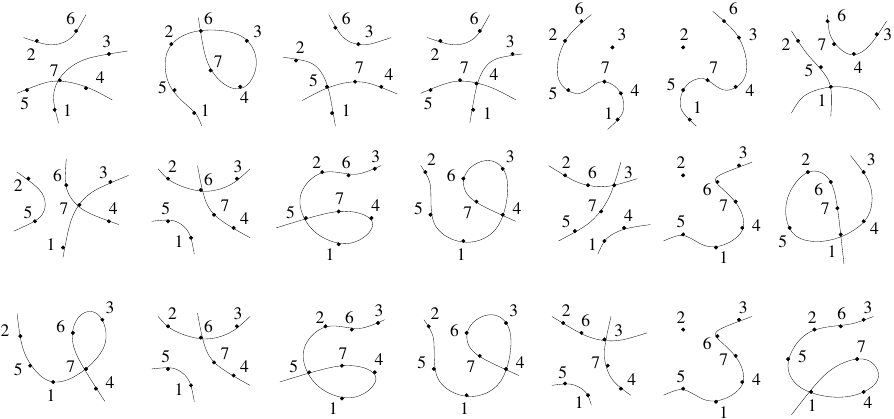}
\caption{ \label{thirdseven} Lists of seven nodal cubics for $R$, $T$, $V$}    
\end{figure}

Three configurations with no six points in convex position may be obtained
taking $7$ in one of the zones $R$, $T$ and $V$ of Figure~\ref{nonconvone}.
The codes of these configurations are shown in Figure~\ref{nonconvtwo}.
We have thus constructed $14$ unordered configurations, to complete the proof of Theorem~2, we need to grant that there exist no others. This will be done in the next section.
The lists of seven nodal cubics corresponding to the three new configurations are
shown in Figure~\ref{thirdseven}. Each cubic may be obtained in several ways perturbing reducible cubics. For example, the cubic with node at $7$ of the first configuration may be obtained from $176 \cup 57423$ moving $6$ to the left of the line $17$, from $274 \cup 56173$ moving $4$ to the top from the line $24$ \dots  
A simple invariant of configurations is obtained counting the numbers of types $\beta, \delta, \gamma$ and $\alpha$ realized by the subcodes $\hat 1, \dots \hat 7$. An encoding for the unordered configurations, using the quadruples $(n_\beta, n_\delta, n_\gamma, n_\alpha)$, is given in Figure~\ref{quadruples}. Note that two unordered configurations have the same quadruple if and only if they are adjacent via a conic-wall. Each configuration is preserved by some subgroup $G$ of $\mathcal{S}_7$, and each $L$-configuration is preserved by a subgroup $G'$. These groups are easily found using the codes. Each point $n$ appears in six subcodes, for each subcode $\hat k$, note the orbit type of $n$ for the action of $G(\hat k)$ (or $G'(\hat k)$). For example, $1$ in $(B, 6)$ realizes $a$ twice and each of the four types $b2$, $d2$, $c1$, $c2$ once. 
An element of $G$ (or $G'$) must map $n$ onto a point $m$ such that: the codes $\hat n$ and $\hat m$ have the same type $\alpha$, $\beta$, $\gamma$ or $\delta$, and $n$, $m$ realize the same distribution of six data (configuration, orbit). Once we have the image $m$ of one point $n$, it is easy to find the images of the other points using the correspondence $\hat n \to \hat m$. This procedure allows to get the groups $G$ and $G'$. The list of groups $G'$ was found originally by Finashin \cite{sf1}, we rediscovered it here with another method. 
In $(\mathbb{R}P^2)^7$ stratified with the lines only, there is a one-to-one correspondence between the topological types for generic elements and the rigid isotopy classes \cite{sf1}. 
More recently, Finashin and Zabun studied the $14$ configurations of seven points 
and proved that the same statement holds for the refined stratification with lines and conics \cite{z}-\cite{fz}.
So, the group $G'$ (resp. $G$) associated to a $L$ (resp. $Q$)-configuration is the {\em monodromy group\/} of $(1, \dots 7) \in (\mathbb{R}P^2)^7$, i.e. the  subgroup of $\mathcal{S}_7$ formed by the permutations of these points realized by $L$ (resp. $Q$)-deformations.





\begin{figure}
\begin{tabular}{ l | l | l | l}
representant & configuration & $G$ & $G'$\\
\hline
$(E, 6)$ & $(7, 0, 0, 0)$ & $\{ id \}$ & $D_7 = \langle (12)(37)(46), (1234567) \rangle$\\
$(D, 6)$ & $(3, 4, 0, 0)_1$ & $\{ id \}$ & $\mathbb{Z}/2 = \langle (34)(25)(16) \rangle$\\
$(G, 6)$ & $(3, 4, 0, 0)_2$ & $\{ id \}$ & $\mathbb{Z}/2 = \langle (14)(23)(57) \rangle$\\
$(C, 6)$ & $(2, 2, 3, 0)_1$ & $\{ id \}$ & $\{ id \}$\\ 
$(C', 6)$ & $(2, 2, 3, 0)_2$ & $\{ id \}$ & $\{ id \}$\\
$(H, 6)$ & $(2, 2, 3, 0)_3$ & $\{ id \}$ & $\{ id \}$\\
$(B, 6)$ & $(1, 2, 2, 2)$ & $\{ id \}$ & $\mathbb{Z}/2 = \langle (16)(25)(34) \rangle$\\
$(A, 6)$ & $(1, 0, 6, 0)$ & $\mathbb{Z}/3 = \langle (135)(246) \rangle$ & $\mathcal{S}_3 = \langle (16)(25)(34), (135)(246) \rangle$\\
$(K, 6)$ & $(1, 6, 0, 0)$ & $\{ id \}$ & $\mathbb{Z}/2 = \langle (14)(36)(25) \rangle$\\
$(I, 6)$  & $(1, 2, 4, 0)$ & $\{ id \}$ & $\mathbb{Z}/2 = \langle (14)(36)(25) \rangle$\\
$(J, 6)$ & $(1, 4, 2, 0)$ & $\{ id \}$ & $\mathbb{Z}/2 = \langle (16)(25)(34) \rangle$\\
$R$ & $(0, 4, 3, 0)$ & $\{ id \}$ & $\{ id \}$\\
$T$ & $(0, 3, 3, 1)$ & $\mathbb{Z}/3 = \langle (137)(456) \rangle$ & $\mathbb{Z}/3 = \langle (137)(456) \rangle$\\
$V$ & $(0, 6, 1, 0)$ & $\mathbb{Z}/3 = \langle (137)(456) \rangle$ & $\mathbb{Z}/3 = \langle (137)(456) \rangle$\\
\end{tabular}
\caption{ \label{quadruples} The $14$ configurations and their monodromy groups}
\end{figure}

\subsection{End of the proof of Theorem~2}
\begin{figure}
\centering
\includegraphics{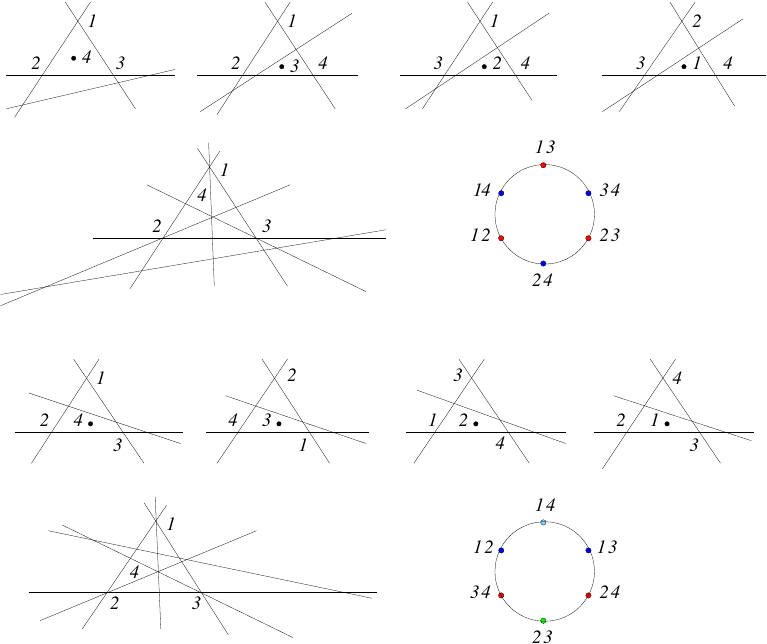}
\caption{ \label{linewalls} The two mutual positions of one line and four points}    
\end{figure}

\begin{figure}
\centering
\includegraphics{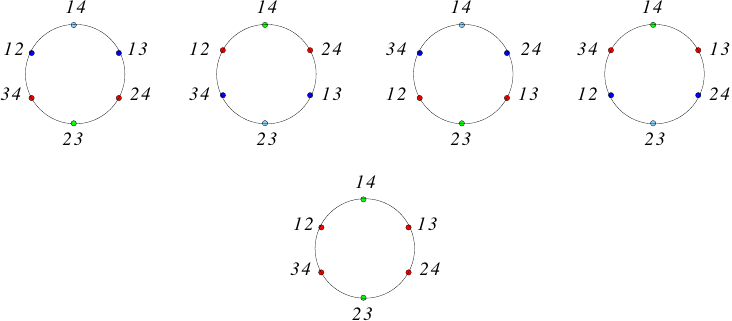}
\caption{ \label{colors} Moving the line}    
\end{figure}

\begin{figure}
\centering
\includegraphics{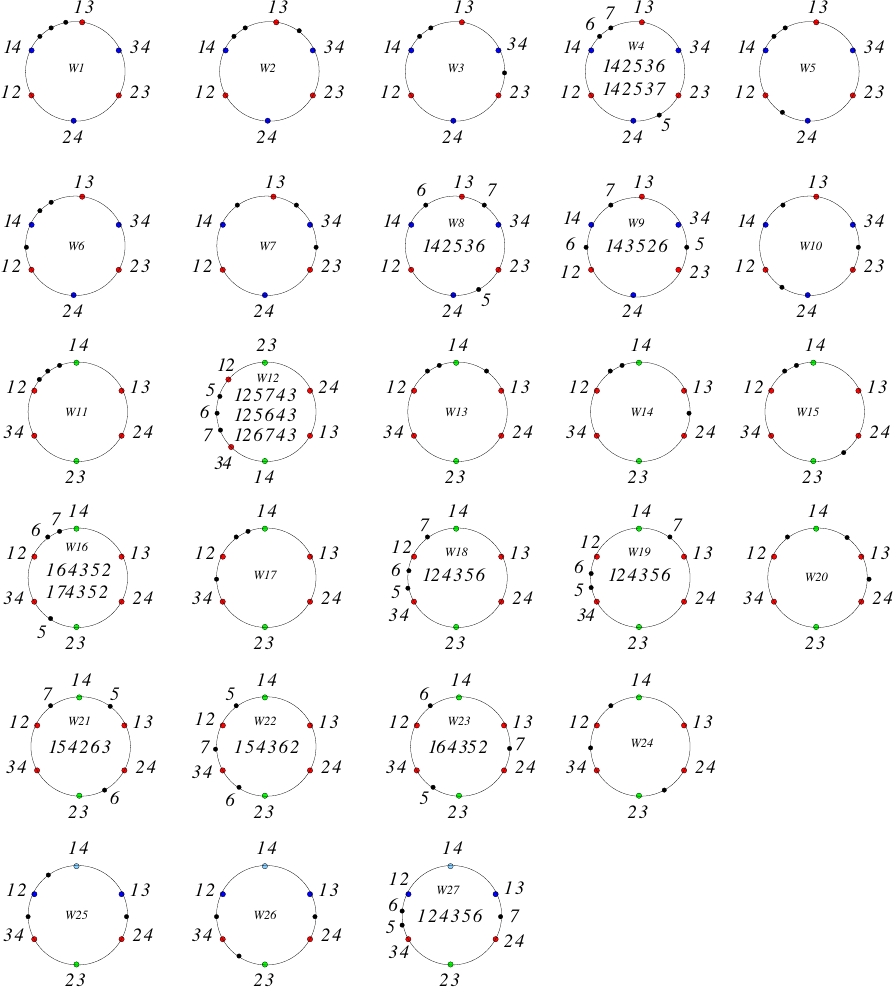}
\caption{ \label{fourtywalls} The $27$ line-walls $W1, \dots W27$, the $14$ line-conic-subwalls}    
\end{figure}

The space $(\mathbb{R}P^2)^7/\mathcal{S}_7$ endowed with the stratification given by the alignment only has $11$ cameras and $27$ walls, see \cite{sf1}, \cite{sf2} (where Finashin considered actually dual arrangements of lines). In our setting, the cameras correspond to the $11$ unordered configurations, described by the quadruples. Let us explain hereafter how to find and encode the walls. 
Consider the possible positions of a line $L$ with respect to four points $1$, $2$, $3$, $4$. Choose three points among $1, \dots 4$, they give rise to four triangles in the plane, we call principal triangle the one containing the fourth point. Add a line $L$ passing through none of the four points, $L$ cuts either three or four principal triangles, see upper and lower part of Figure~\ref{linewalls}. 
The corresponding symmetry groups are $\mathcal{S}_3 = \langle (12), (123) \rangle$ and $\mathbb{Z}/2 = \{ $id$, (23) \}$.
The four points give rise to six lines, let us look at the six intersection points of these lines with $L$. The cyclic ordering of these intersections on $L$ allows to recover the mutual position of $L$ with $1, 2, 3, 4$. This information will be encoded using a circle with six marked points, two points have the same color if they are in the same orbit for the action of the group. We need two colors (red and blue) in the first case, and four colors (pale blue, dark blue, red and green) in the second case. To get the walls in $(\mathbb{R}P^2)^7/\mathcal{S}_7$, distribute three unmarked points (colored black) in all possible ways on the two circles from Figure~\ref{linewalls}, that stay for first marked. For the second circle, if an interval between a blue point, say $24$, and a red point, say $13$, contains no black point, we may move the line $L$ (or equivalently the point $4$) until $L$ crosses the intersection $24 \cap 13$ (doing so, we won't leave the wall). On the circle, the positions of $13$ and $24$ are swapped, and the colorings of all six points change, see second circle in Figure~\ref{colors}. The third and fourth circle show the colorings obtained with the last two possible positions of $34$, $12$, $13$, $24$.
In other words, the group $G$ acting on the four points is $D_4 =
\{ $id$, (23), (14), (12)(34), (13)(24), (14)(23), (1243), (1342) \}$.
We can thus reduce the number of colors to two (red and green), see Figure~\ref{colors}, bottom circle. 
(We may see the pairs $2, 3$ and $1, 4$ as opposite edges of a square, the red points on the circle are the vertices of the square, the green points are the intersections of the pairs of opposite lines supporting the edges).
Remove now the markings on the circles, one gets in total $27$ line-walls $W1, \dots W27$, splitting in three groups according to the distribution of colors for the six intersections with the line $L$, see the $27$ circles in Figure~\ref{fourtywalls} (where the markings should be ignored, they will be used in the next section).
The combinatorial types we have defined here are rigid isotopy invariants
for the elements of $(\mathbb{R}P^2)^7$ with three aligned points. Conversely, 
consider two $7$-uples of points $(P_1, \dots, P_7)$ and $(P'_1, \dots P'_7)$, realizing the same combinatorial type.
A projective transformation $\sigma$ maps $(L'_{12}, L'_{13}, L'_{23}, L')$ onto $(L_{12}, L_{13}, L_{23}, L)$. One has $\sigma(P'_i) = P_i$, $i = 1, 2, 3$, and
$\sigma(P'_4)$ lies in the zone of $\mathbb{R}P^2\setminus (L_{12} \cup L_{13} \cup L_{23} \cup L)$ containing $P_4$.  A further deformation moves then $\sigma(P'_4)$ onto $P_4$ (inducing possibly a change as shown in Figure~\ref{colors}), and shifts the images of $P'_5$, $P'_6$, $P'_7$ along $L$ onto $P_5$, $P_6$, $P_7$.

\begin{proposition}
For each of the $14$ configurations $(E, 6), \dots V$, the unordered configurations realized by the adjacent configurations are shown in Figures~\ref{adj1}-\ref{adj4}, along with the corresponding walls. 
\end{proposition}

{\em Proof:\/} 
For each configuration $(E, 6), \dots V$, we find out the list of adjacent line-walls in form of triples: the rule from Figures~\ref{movebeta}-\ref{movealpha} allows to get all of the adjacent triples for each subcode $\hat 1, \dots \hat 7$, we get thus seven lists of triples. The adjacent triples $lnm$ for the configuration are those appearing four times in total (in all lists except for $\hat l$, $\hat n$, $\hat m$). Draw the code of the configuration obtained after the crossing. Note that we need the markings only in the case that the configuration obtained is of type $(3, 4, 0, 0)$ or $(2, 2, 3, 0)$.
Each crossing induces the change of four subcodes. To determine the type of the line-wall among $W1, \dots W27$, we need to determine the cyclic ordering of nine points on the line $L$. Write, for each of the four subcodes, the cyclic ordering of the six relevant points on $L$ using the rule of Figures~\ref{movebeta}-\ref{movealpha}. We get thus four circles marked with six points each. The partial cyclic orderings given by these four circles allow to recover the cyclic ordering of the nine points without ambiguity, except for pairs of adjacent red points. For example, start from $(D, 6)$ and take the crossing $456$, the partial cyclic orderings are: $4, 5, 6, 27, 37, 23$  (change of $\hat 1$), $4, 5, 17, 6, 37, 13$ ($\hat 2$), $4, 5, 17, 6, 27, 12$ ($\hat 3$), $4, 5, 6, 12, 13, 23$ ($\hat 7$).  The resulting cyclic ordering is
$4, 5, 17, 6, 27, $\{37, 12\}$, 13, 23$, the wall is $W18$.
The informations obtained are gathered in Figures~\ref{adj1}-\ref{adj3}. 
By Proposition~2, there are eleven unordered configurations  with six coconic points, they have representants $A, \dots K$. Each ordered configuration $(X, 6)$ is adjacent to $(X, 6')$ via the conic-wall $X$. The equivalences between ordered configuration are: $(X, 6) \simeq (X, 6')$ for $X \in \{ A, B, I, J, K \}$; $(D, 6) \simeq (D, 6') \simeq (G, 6')$; $(E, 6) \simeq (E, 6') \simeq (F, 6) \simeq (F, 6')$; and $(C, 6) \simeq (H, 6')$. We get thus the informations of Figure-\ref{adj4}, Proposition~3 is proved. It follows from this Proposition that there exist no other unordered configurations than the $14$, this finishes the proof of Theorem~2. $\Box$

\begin{figure}
\centering
\includegraphics{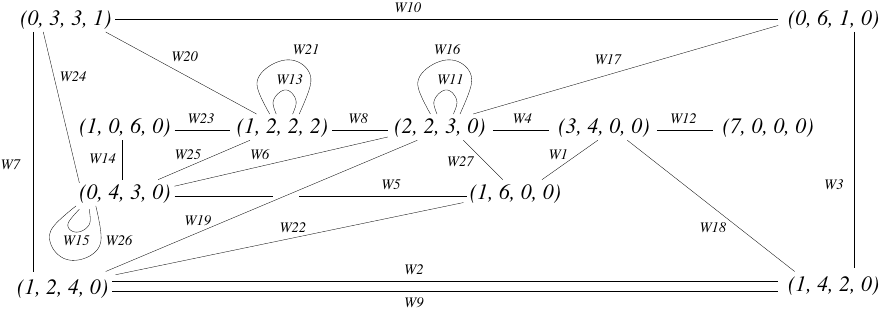}
\caption{ \label{adjacency} Adjacency graph for seven unordered points, stratification by lines only}    
\end{figure}

The adjacency graph of the space  $(\mathbb{R}P^2)^7/\mathcal{S}_7$, endowed with the stratification given by the alignment only, is constructed in \cite{sf1}, \cite{sf2}. With help of Figures~\ref{adj1}-\ref{adj3}, we recover this graph, see Figure~\ref{adjacency}. For each arrangement, Finashin obtained the adjacent line-walls as orbits of triangles under the action of its group of symmetries.
In our setting, the arrangements of lines and triangles are replaced by  configurations of points and triples. 
We choose representants of the eleven $L$-configurations, for example by removing $(C'(6)$, $(G, 6)$ and $(H, 6)$ from our canonical list of fourteen $Q$-configurations. Encode these representants now as $L$-configurations, removing the colors black and white from the dots on the circles representing the $\beta$-codes. Now, each $\beta$-code has six adjacent triples instead of three, namely the triples of three consecutive points on the circle. (For the $\delta$, $\gamma$ and $\alpha$-codes, the sets of adjacent triples remain unchanged.) We may use the same rule as at the beginning of this section to get the adjacent triples $lnm$ for the eleven $L$-configurations. The line-walls adjacent to each $L$-configuration correspond to the orbits of these triples under the action of its monodromy group $G'$, see Figure~\ref{quadruples}. For example $(A, 6)$ has two adjacent line-walls: $\{ 367, 257, 147 \}$ ($W23$) and $\{ 234, 456, 126, 345, 123, 156 \}$ ($W14$). There is one exception: the configuration $(0, 4, 3, 0)$, represented by $R$, has monodromy group $\{ id \}$, but the adjacent wall $W26$ is represented by two triples, $125$ and $357$. As a matter of fact, $W26$ is a two-sided inner wall. The other five inner walls $W15$, $W13$, $W21$, $W11$ and $W16$ are one-sided \cite{sf1}.

\begin{figure}
\begin{tabular}{ c c c c c }
$(E, 6)$ & $(7, 0, 0, 0)$ & $167$ & $W12$ & $(3, 4, 0, 0)_1$\\
 & & $456$ & $W12$ & $(3, 4, 0, 0)_2$\\
 & & & & \\
$(D, 6)$ & $(3, 4, 0, 0)_1$ & $267$ & $W4$ & $(2, 2, 3, 0)_1$\\
 & & $157$ & $W4$ & $(2, 2, 3, 0)_2$\\
 & & $456$ & $W18$ & $(1, 4, 2, 0)$\\
 & & $167$ & $W12$ & $(7, 0, 0, 0)$\\
 & & & & \\
$(C, 6)$ & $(2, 2, 3, 0)_1$ & $267$ & $W4$ & $(3, 4, 0, 0)_1$\\
 & & $157$ & $W8$ & $(1, 2, 2, 2)$\\
 & & $367$ & $W16$ & $(2, 2, 3, 0)_1$\\
 & & $126$ & $W27$ & $(1, 6, 0, 0)$\\
 & & $456$ & $W19$ & $(1, 2, 4, 0)$\\
 & & & & \\
$(B, 6)$ & $(1, 2, 2, 2)$ & $157$ & $W8$ & $(2, 2, 3, 0)_1$\\
 & & $267$ & $W8$ & $(2, 2, 3, 0)_2$\\
 & & $257$ & $W23$ & $(1, 0, 6, 0)$\\
 & & $367$ & $W21$ & $(1, 2, 2, 2)$\\
 & & $234$ & $W13$ & $(1, 2, 2, 2)$\\
 & & $456$ & $W20$ & $(0, 3, 3, 1)$\\
 & & $147$ & $W21$ & $(1, 2, 2, 2)$\\
 & & $126$ & $W25$ & $(0, 4, 3, 0)$\\
 & & & & \\
$(A, 6)$ & $(1, 0, 6, 0)$ & $257$ & $W23$ & $(1, 2, 2, 2)$\\
 & & $367$ & $W23$ & $(1, 2, 2, 2)$\\ 
 & & $234$ & $W14$ & $(0, 4, 3, 0)$\\
 & & $456$ & $W14$ & $(0, 4, 3, 0)$\\
 & & $147$ & $W23$ & $(1, 2, 2, 2)$\\
 & & $126$ & $W14$ & $(0, 4, 3, 0)$\\
 & & & & \\
$C'(6)$ & $(2, 2, 3, 0)_2$ & $157$ & $W4$ & $(3, 4, 0, 0)_1$\\
 & & $267$ & $W8$ & $(1, 2, 2, 2)$\\
 & & $234$ & $W11$ & $(2, 2, 3, 0)_2$\\
 & & $456$ & $W17$ & $(0, 6, 1, 0)$\\
 & & $147$ & $W16$ & $(2, 2, 3, 0)_2$\\
\end{tabular}
\caption{Adjacencies via line-walls \label{adj1} }
\end{figure}

\begin{figure}
\begin{tabular}{ c c c c c }
$(G, 6)$ & $(3, 4, 0, 0)_2$ & $234$ & $W1$ & $(1, 6, 0, 0)$\\
 & & $567$ & $W12$ & $(7, 0, 0, 0)$\\
 & & $467$ & $W4$ & $(2, 2, 3, 0)_3$\\
 & & $127$ & $W18$ & $(1, 4, 2, 0)$\\
 & & & & \\
$(H, 6)$ & $(2, 2, 3, 0)_3$ & $467$ & $W4$ & $(3, 4, 0, 0)_2$\\
 & & $457$ & $W27$ & $(1, 6, 0, 0)$\\
 & & $234$ & $W6$ & $(0, 4, 3, 0)$\\
 & & $367$ & $W16$ & $(2, 2, 3, 0)_3$\\
 & & $127$ & $W19$ & $(1, 2, 4, 0)$\\
& & & & \\
$(K, 6)$ & $(1, 6, 0, 0)$ & $456$ & $W1$ & $(3, 4, 0, 0)_2$\\
 & & $234$ & $W5$ & $(0, 4, 3, 0)$\\
 & & $367$ & $W22$ & $(1, 2, 4, 0)$\\
 & & $457$ & $W27$ & $(2, 2, 3, 0)_3$\\
 & & $127$ & $W27$ & $(2, 2, 3, 0)_1$\\
& & & & \\
$(I, 6)$ & $(1, 2, 4, 0)$ & $457$ & $W19$ & $(2, 2, 3, 0)_1$\\
 & & $234$ & $W7$ & $(0, 3, 3, 1)$\\
 & & $467$ & $W9$ & $(1, 4, 2, 0)$\\
 & & $367$ & $W22$ & $(1, 6, 0, 0)$\\
 & & $137$ & $W9$ & $(1, 4, 2, 0)$\\
 & & $126$ & $W2$ & $(1, 4, 2, 0)$\\
 & & $127$ & $W19$ & $(2, 2, 3, 0)_3$\\
 & & & & \\
$(J, 6)$ & $(1, 4, 2, 0)$ & $467$ & $W9$ & $(1, 2, 4, 0)$\\
 & & $234$ & $W2$ & $(1, 2, 4, 0)$\\
 & & $567$ & $W18$ & $(3, 4, 0, 0)_1$\\
 & & $137$ & $W9$ & $(1, 2, 4, 0)$\\
 & & $126$ & $W3$ & $(0, 6, 1, 0)$\\
 & & $127$ & $W18$ & $(3, 4, 0, 0)_2$\\ 
\end{tabular}
\caption{Adjacencies via line-walls, continued \label{adj2} }
\end{figure}

\begin{figure}
\begin{tabular}{ c c c c c }
$R$ & $(0, 4, 3, 0)$ & $236$ & $W25$ & $(1, 2, 2, 2)$\\
 & & $457$ & $W6$ & $(2, 2, 3, 0)_3$\\
 & & $247$ & $W14$ & $(1, 0, 6, 0)$\\
 & & $357$ & $W26$ & $(0, 4, 3, 0)$\\
 & & $167$ & $W15$ & $(0, 4, 3, 0)$\\
 & & $156$ & $W5$ & $(1, 6, 0, 0)$\\
 & & $134$ & $W24$ & $(0, 3, 3, 1)$\\
 & & $125$ & $W26$ & $(0, 4, 3, 0)$\\
 & & & & \\
$T$ & $(0, 3, 3, 1)$ & $247$ & $W20$ & $(1, 2, 2, 2)$\\
 & & $356$ & $W7$ & $(1, 2, 4, 0)$\\
 & & $236$ & $W20$ & $(1, 2, 2, 2)$\\
 & & $357$ & $W24$ & $(0, 4, 3, 0)$\\
 & & $467$ & $W7$ & $(1, 2, 4, 0)$\\
 & & $145$ & $W7$ & $(1, 2, 4, 0)$\\
 & & $137$ & $W10$ & $(0, 6, 1, 0)$\\
 & & $167$ & $W24$ & $(0, 4, 3, 0)$\\
 & & $125$ & $W20$ & $(1, 2, 2, 2)$\\
 & & $134$ & $W24$ & $(0, 4, 3, 0)$\\
 & & & & \\
$V$ & $(0, 6, 1, 0)$ & $247$ & $W17$ & $(2, 2, 3, 0)_2$\\
 & & $356$ & $W3$ & $(1, 4, 2, 0)$\\
 & & $236$ & $W17$ & $(2, 2, 3, 0)_2$\\
 & & $467$ & $W3$ & $(1, 4, 2, 0)$\\
 & & $145$ & $W3$ & $(1, 4, 2, 0)$\\
 & & $137$ & $W10$ & $(0, 3, 3, 1)$\\
 & & $125$ & $W17$ & $(2, 2, 3, 0)_2$\\
\end{tabular}
\caption{Adjacencies via line-walls, end \label{adj3} }
\end{figure}

\begin{figure}
\begin{tabular}{ c c c c }
$(E, 6)$ & $(7, 0, 0, 0)$ & $E$ & $(7, 0, 0, 0)$\\
 & & $F$ & $(7, 0, 0, 0)$\\
$(D, 6)$ & $(3, 4, 0, 0)_1$ & $D$ & $(3, 4, 0, 0)_1$\\
 & & $G$ & $(3, 4, 0, 0)_2$\\
$(C, 6)$ & $(2, 2, 3, 0)_1$ & $C$ & $(2, 2, 3, 0)_2$\\
 & & $H$ & $(2, 2, 3, 0)_3$\\
$(B, 6)$ & $(1, 2, 2, 2)$ & $B$ & $(1, 2, 2, 2)$\\
$(A, 6)$ & $(1, 0, 6, 0)$ & $A$ & $(1, 0, 6, 0)$\\
$(C', 6)$ & $(2, 2, 3, 0)_2$ & $C$ & $(2, 2, 3, 0)_1$\\
$(G, 6)$ & $(3, 4, 0, 0)_2$ & $G$ & $(3, 4, 0, 0)_1$\\
$(H, 6)$ & $(2, 2, 3, 0)_3$ & $H$ & $(2, 2, 3, 0)_1$\\
$(K, 6)$ & $(1, 6, 0, 0)$ & $K$ & $(1, 6, 0, 0)$\\
$(I, 6)$ & $(1, 2, 4, 0)$ & $I$ & $(1, 2, 4, 0)$\\
$(J, 6)$ & $(1, 4, 2, 0)$ & $J$ & $(1, 4, 2, 0)$\\
\end{tabular}
\caption{Adjacencies via conic-walls \label{adj4} }
\end{figure}

\subsection{Configurations with three aligned points}
\begin{proposition}
Seven points in $\mathbb{R}P^2$, with three of them aligned, but otherwise generic,
may realize $38$ different unordered configurations. 
\end{proposition}
{\em Proof:\/} 
Let us call {\em refined line-walls\/} the configurations with three aligned points, otherwise generic. To find them, we need first to determine, for each of the $27$ line-walls, all sets of six points that may be coconic. For each line-wall with some admissible conic, we mark the black points with names $5, 6, 7$, and indicate the conic(s) inside of the circle, see  Figure~\ref{fourtywalls}. Let us explain how to spot these conics.  A non-generic configuration is {\em line-conic\/} if it has a conic through six points and a line through three points, one of them not on the conic. The $11$ unordered configurations with six coconic points are represented by $A, \dots K$, with conic $123456$. 
Figure~\ref{triples} displays all triples of points with $7$ that may become aligned for $A, \dots K$. Most of these triples are directly visible in Figure~\ref{hexagon}, but the safest way to find them without forgetting one is the following: write the list of all admissible triples for each subconfiguration $\hat 1, \dots \hat 6$ using the rules explicited in Figures~\ref{movebeta}-\ref{movealpha}, the relevant triples $nm7$ are those appearing four times in total (in all lists but $\hat n, \hat m$). Each of the line-conic configurations obtained may be denoted by a letter followed by a triple, we find in total $36$ combinations letter-triple.  Note that a line-conic configuration may be encoded with two letter-triples, for example, $A357 = B357$, see Figure~\ref{hexagon}. Note also that some are mapped onto others by elements of $\mathcal{S}_7$: the cyclic permutation $(1 6 5 4 3 2)$ sends $A367$ onto $A257$ and $A257$ onto $A147$; the symmetry $s = (1 6)(2 5)(3 4)$ maps $Xij7$ onto $Xs(i)s(j)7$ for $X \in \{ A, B, D, E, F, J \}$. The symmetry $(1 4)(2 5)(3 6)$ swaps the elements of the pairs ($I457$, $H127$),  ($K127$, $H457$) and ($I137$, $I467$).
We find in total $14$ different line-conic unordered configurations, for each of them, choose a letter-triple as representant. See then Figure~\ref{hexagon} to get the cyclic ordering with which the line meets the three aligned points, and the six lines determined by the other four points. 
The left column of Figure~\ref{conicwalls} displays the symmetries mapping the $14$ letter-triples onto the conic-subwalls in Figure~\ref{fourtywalls}. The right column displays the lists of equivalent letter-triples.


\begin{figure}
\begin{tabular}{ c | l }
$A$ & $367$, $257$, $147$\\
$B$ & $367$, $267$, $257$, $147$, $157$\\
$C$ & $367$, $267$, $157$\\
$D$ & $157$, $267$, $167$\\
$E$ & $167$\\
$F$ & $127$, $567$\\
$G$ & $467$, $567$, $127$\\
$H$ & $367$, $457$, $467$, $127$\\
$I$ & $367$, $457$, $467$, $137$, $127$\\
$J$ & $467$, $567$, $137$, $127$\\
$K$ & $457$, $367$, $127$\\
\end{tabular}
\caption{ \label{triples} Admissible triples for the $11$ configurations with six coconic points}
\end{figure}

\begin{figure}
\begin{tabular}{ c  c  c | l }
$C267$ & $(1 3 2 5)$ & $W4, 142536$ & $D267$, $D157$\\
$G467$ & $(2 4 5 3)(6 7)$ & $W4, 142537$ & $H467$\\
$B157$ & $(1 6 3 4 2)$ & $W8, 142536$ & $B267$, $C157$\\
$I467$ & $(2 4 5)$ & $W9, 143526$ & $J467$, $J137$, $I137$\\
$D167$ & $(1 7 6 5 2 4)$ & $W12, 125743$ & $E167$\\
 & $(1 5 4 3)(6 7)$ &  & \\
$F127$ & $(1 6 2 7 5)(3 4)$ & $W12, 126743$ & $F567$, $G567$\\
 & $(1 6 4)(2 5 3)$ & $W12, 125643$ & \\
$C367$ & $(3 5 4)$ & $W16, 164352$ & \\
 & $(1 4 2 3 5)$ &  & \\
$H367$ & $(3 5 4)(6 7)$ & $W16, 174352$ & \\
 & $(1 4 2 3 5)(6 7)$ &  & \\
$G127$ & $(1 6)(2 5)$ & $W18, 124356$ & $J127$, $J567$ \\
$H127$ & $(1 6)(2 5)$ & $W19, 124356$ & $I127$, $I457$\\
$B367$ & $(2 3 6 5 4)$ & $W21, 154263$ & $B147$\\
 & $(1 4 3 6 5)$ & & \\
$I367$ & $(1 3 5 2 4)$ & $W22, 154362$ & $K367$\\
 & $(3 6 5 4)$ & & \\
$A367$ & $(1 3 6 5 2 4)$ & $W23, 164352$ & $A257$, $A147$, $B257$\\
 & $(1 4 2 3 5)$ & & \\
$H457$ & $(1 2 4 5 6)$ & $W27, 124356$ & $K127$, $K457$\\
\end{tabular}
\caption{ \label{conicwalls} The $14$ line-conic-walls}
\end{figure}

\begin{figure}
\centering
\includegraphics{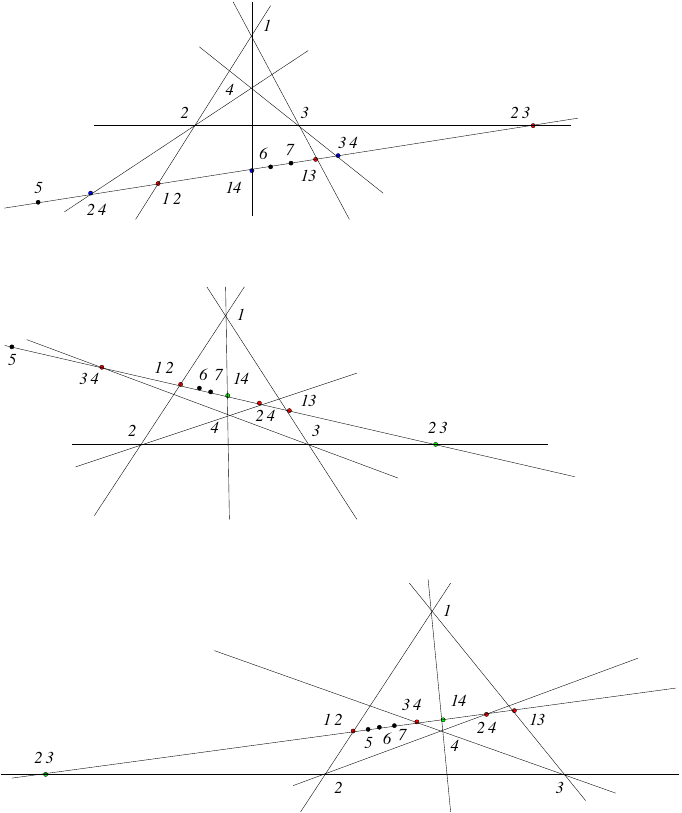}
\caption{ \label{multwalls} Line-walls $W4$, $W16$ and $W12$}    
\end{figure}

Let us look at Figure~\ref{fourtywalls}. There are actually $15$ conics in total, distributed in $11$ line-walls. But the two conic-subwalls $W12, 125743$, $W12, 126743$ are swapped by the symmetry $(1 3)(2 4)(7 5)$, hence equivalent. (Note also that this symmetry maps $W12, 125743$ onto itself). 
There are eight line-walls with one unique conic subwall. Each of them gives rise to two refined line-walls, obtained moving one of the six points away from the conic either to the inside or to the outside, so as to realize a $\beta$ configuration. We shall encode it as follows. Say $123456$ lie in convex position and each of the points $1, 3, 5$ is exterior to the conic through the other five, we write shortly $(123456, 135)$.
The symmetry $(1 3)(2 4)(5 6)$ maps $W22, 154362$ onto itself, and swaps the refined-subwalls $W22, (154362, 235)$ with $W22, (154362, 146)$.
This symmetry maps also $W23, 164352$ onto itself, and swaps the refined subwalls $W23, (164352, 236)$ with $W23, (164352, 145)$.
The three line-walls $W4$, $W12$ and $W16$ that have several conic-subwalls are represented in Figure~\ref{multwalls}. 
Let $n$ be a point and $C_2$ be a conic, we write $n < C_2$ ($n > C_2$) if $n$ lies inside (outside) of $C_2$.
For $W4$, one has: $6, 7 < 14253$ or $6 > 14253$ and $7 < 14253$, or $6, 7 > 14253$.
For $W16$, one has: $6, 7 < 14352$ or $6 > 14352$ and $7 < 14352$, or $6, 7 > 
14352$.
Each of the two line-walls $W4$ and $W16$ has three refined line-walls, see the encoding with $\beta$-codes in Figure~\ref{reflinewalls}.
The wall $W12$ has four refined line-walls, corresponding to the cases: $6, 7 < 12543$ and $6 < 12743$; $6, 7 < 12543$ and $6 > 12743$; $6 < 12543$ and $7 > 12543$; $6, 7 > 12543$. The triples of $\beta$ codes $\hat 5$, $\hat 6$ and $\hat 7$ for these four cases are:

$\hat 5 = 126743, 273$, $\hat 6 = 125743, 154$, $\hat 7 = 125643, 154$,  ($125743$, $126743$) 

$\hat 5 = 126743, 164$, $\hat 6 = 125743, 154$, $\hat 7 = 125643, 154$,  ($126743$) 

$\hat 5 = 126743, 273$, $\hat 6 = 125743, 273$, $\hat 7 = 125643, 154$, ($125643$, $125743$),

$\hat 5 = 126743, 273$, $\hat 6 = 125743, 273$, $\hat 7 = 125643, 263$, ($125643$)\

\begin{figure}
\centering
\includegraphics{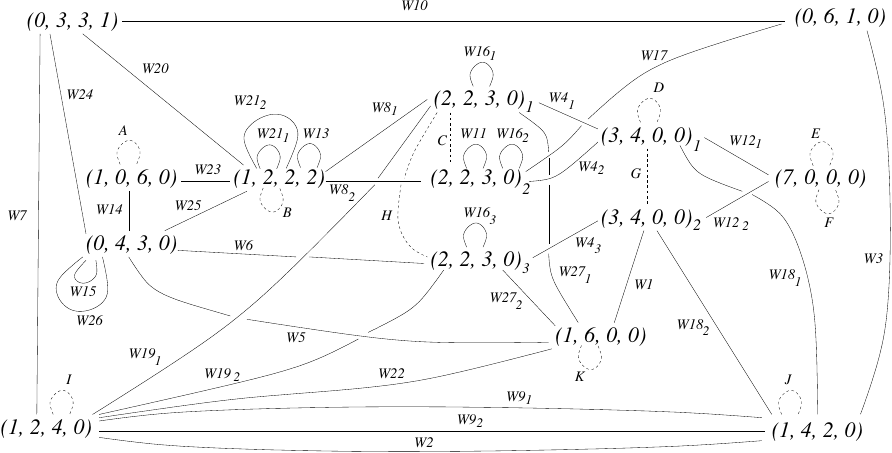}
\caption{ \label{adjaref} Adjacency graph for seven unordered points, stratification by lines and conics}    
\end{figure}

(We have indicated the adjacent conic-walls for each of them.)
The symmetry $(1 3)(2 4)(5 7)$ swaps the refined line-walls pairwise: the first one with the third, and the second one with the fourth.
So $W12$ gives rise to only two unordered refined line-walls.
To characterize them, let us consider the middle point in the group of three aligned points ($6$ in Figure~\ref{fourtywalls}). Either $6$ is interior for one only of the codes $\hat 5$, $\hat 7$, or $6$ is interior for both codes $\hat 5$ and $\hat 7$.
The $11$ line-walls with conic-subwalls give rise together to $22$ particular unordered refined line-walls. A representant for each of them is obtained adding some $\beta$-codes after the name $W_N$ of the wall. 
For $N = 12$, we need three $\beta$-codes, for $N = 4, 16$ we need two, for the other values of $N$, we need only one, see Figure~\ref{reflinewalls}.
We have yet to find out, for each refined line-wall, the adjacent pair of unordered configurations. Appropriate symmetries mapping our $22$ representants onto walls between configurations of types $(X, 6)$ or $(X, 6')$ allow to achieve this, see Figures~\ref{reflinewalls}~\ref{refconf}.
The rigid isotopy classes for the elements of $(\mathbb{R}P^2)^7$ with three aligned points correspond one-to-one to the refined combinatorial types. To prove this, the same argument as for the stratification by lines only applies, combined with the fact that projective transformations preserve the mutual positions of points and conics.
The total number of refined line-walls is $22 + (27 - 11) = 38$.
The adjacency graph for the space of seven unordered points, stratified by lines and conics, is displayed in  Figure~\ref{adjaref}.
Here again, the walls adjacent to a configuration correspond to the orbits of the adjacent triples and eventual adajacent cyclically ordered sextuples under the action of its monodromy group $G$. 
There is one single exception: the two-sided inner line-wall $W26$, see end of section 3.2. All of the other inner walls are one-sided.
$\Box$ 


Rational cubics and pencils of rational cubics were used in \cite{fi1}-\cite{fi3} to study the topology of real algebraic $M$-curves of degree $9$. The combinatorics of generic pencils of cubics was studied in \cite{fi4}, for the particular case where eight of the base points lie in convex position.

\begin{figure}
\begin{tabular}{ c  c  c  c c  c }  
$W4_1$ & $(142536, 456)$  & $(1 5 2 3)$ & $(C, 6)$ & $267$ & $(D, 6)$\\
            & $(142537, 123)$  &                  &              &           &       \\
$W4_2$ & $(142536, 123)$  & $(1 5 2 3)$ & $(C, 6')$ & $267$ & $(D, 6')$\\
            & $(142537, 123)$  &                  &              &           &       \\
$W4_3$ & $(142536, 456)$  & $(2 3 5 4)(6 7)$ & $(H, 6)$ & $467$ & $(G, 6)$\\
            & $(142537, 457)$ &                  &              &           &       \\
$W8_1$ & $(142536, 123)$ & $(1 2 4 3 6)$ & $(B, 6)$ & $157$ & $(C, 6)$\\
$W8_2$ & $(142536, 456)$ & $(1 2 4 3 6)$ & $(B, 6')$ & $157$ & $(C, 6')$\\
$W9_1$ & $(143526, 456)$ & $(2 5 4)$ & $(I, 6)$ & $467$ & $(J, 6)$\\
$W9_2$ & $(143526, 132)$ & $(2 5 4)$ & $(I, 6')$ & $467$ & $(J, 6')$\\
$W12_1$ & $(126743, 273)$  & $(1 3 4 5)(6 7)$ & $(E, 6)$ & $167$ & $(D, 6)$\\
            &  $(125743, 154)$ &                  &              &           &       \\
            &  $(125643, 154)$ &                  &              &           &       \\
$W12_2$ & $(126743, 164)$   & $(1 3 2 4)$ & $(F, 6)$ & $567$ & $(G, 6)$\\
            & $(125743, 154)$  &                  &              &           &       \\
            & $(125643, 154)$ &                  &              &           &       \\
$W16_1$ & $(164352, 236)$  & $(3 4 5)$ & $(C, 6)$ & $367$ & $(1 5)(2 4)(C, 6)$\\
            & $(174352, 145)$ &                  &              &           &       \\
$W16_2$ & $(164352, 145)$  & $(3 4 5)$ & $(C, 6')$ & $367$ & $(1 5)(2 4)(C, 6')$\\
            & $(174352, 145)$ &                  &              &           &       \\
$W16_3$ & $(164352, 236)$  & $(6 7)(3 4 5)$ & $(H, 6)$ & $367$ & $(1 5)(2 4)(H, 6)$\\
            & $(174352, 237)$ &                  &              &           &       \\
$W18_1$ & $(124356, 236)$ & $(1 6)(2 5)$ & $(G, 6')$ & $127$ & $(J, 6')$\\
$W18_2$ & $(124356, 145)$ & $(1 6)(2 5)$ & $(G, 6)$ & $127$ & $(J, 6)$\\
$W19_1$ & $(124356, 236)$ & $(1 6)(2 5)$ & $(H, 6')$ & $127$ & $(I, 6')$\\
$W19_2$ & $(124356, 145)$ & $(1 6)(2 5)$ & $(H, 6)$ & $127$ & $(I, 6)$\\
$W21_1$ & $(154263, 235)$ & $(2 4 5 6 3)$ & $(B, 6)$ & $367$ & $(1 5)(2 4)(B, 6)$\\
$W21_2$ & $(154263, 146)$ & $(2 4 5 6 3)$ & $(B, 6')$ & $367$ & $(1 5)(2 4)(B, 6')$\\
$W22$ & $(154362, 146)$ & $(1 4 2 5 3)$ & $(I, 6)$ & $367$ & $(K, 6)$\\
$W23$ & $(164352, 236)$ & $(1 5 3 2 4)$ & $(A, 6)$ & $367$ & $(B, 6)$\\
$W27_1$ & $(124356, 236)$ & $(1 6 5 4 2)$ & $(H, 6')$ & $457$ & $(K, 6')$\\
$W27_2$ & $(124356, 145)$ & $(1 6 5 4 2)$ & $(H, 6)$ & $457$ & $(K, 6)$\\
\end{tabular}
\caption{ \label{reflinewalls} Adjacencies via the $22$ particular refined line-walls}
\end{figure}

\begin{figure}
\begin{tabular}{ c  c  c }
$(3, 4, 0, 0)_1$ & $W4_1$ & $(2, 2, 3, 0)_1$\\
                        & $W4_2$ & $(2, 2, 3, 0)_2$\\
$(3, 4, 0, 0)_2$ & $W4_3$ & $(2, 2, 3, 0)_3$\\
$(1, 2, 2, 2)$ & $W8_1$ & $(2, 2, 3, 0)_1$\\
                    & $W8_2$ & $(2, 2, 3, 0)_2$\\
$(1, 2, 4, 0)$ & $W9_1$ & $(1, 4, 2, 0)$\\
                    & $W9_2$ & $(1, 4, 2, 0)$\\
$(7, 0, 0, 0)$ & $W12_1$ & $(3, 4, 0, 0)_1$\\
                    & $W12_2$ & $(3, 4, 0, 0)_2$\\
$(2, 2, 3, 0)_1$ & $W16_1$ & $(2, 2, 3, 0)_1$\\
$(2, 2, 3, 0)_2$ & $W16_2$ & $(2, 2, 3, 0)_2$\\
$(2, 2, 3, 0)_3$ & $W16_3$ & $(2, 2, 3, 0)_3$\\
$(1, 4, 2, 0)$ & $W18_1$ & $(3, 4, 0, 0)_1$\\
                    & $W18_2$ & $(3, 4, 0, 0)_2$\\
$(1, 2, 4, 0)$ & $W19_1$ & $(2, 2, 3, 0)_1$\\
                    & $W19_2$ & $(2, 2, 3, 0)_3$\\
$(1, 2, 2, 2)$ & $W21_1$ & $(1, 2, 2, 2)$\\
                    & $W21_2$ & $(1, 2, 2, 2)$\\
$(1, 6, 0, 0)$ & $W22$ & $(1, 2, 4, 0)$\\
$(1, 0, 6, 0)$ & $W23$ & $(1, 2, 2, 2)$\\
$(1, 6, 0, 0)$ & $W27_1$ & $(2, 2, 3, 0)_1$\\
                    & $W27_2$ & $(2, 2, 3, 0)_3$\\
\end{tabular}
\caption{ \label{refconf} Adjacencies between unordered configurations}
\end{figure}

\vspace{2cm}
severine.fiedler@live.fr

\end{document}